\documentclass[12pt]{article}
\usepackage[utf8]{inputenc}
\usepackage{graphicx}
\usepackage{amsmath, amssymb,theorem,verbatim}
\usepackage{hyperref}
\usepackage{textcomp} 
\usepackage{setspace}
\usepackage{mathtools}
\usepackage[round]{natbib}
\usepackage{multirow} 
\usepackage{authblk}
\usepackage[usenames,dvipsnames]{xcolor}

\numberwithin{equation}{section}

\newtheorem{prop}{Proposition}[section] 
\newtheorem{theorem}{Theorem}[section]
\newtheorem{lemma}{Lemma}[section] 
\newtheorem{assumption}{Assumption}[section] 
 
\newtheorem{corollary}{Corollary}[section] 

\newtheorem{remark}{Remark}[section]

\newcommand{\C}{\mathbb{C}}

\newcommand{\N}{\mathbb{N}}
\newcommand{\R}{\mathbb{R}}

\newcommand{\Z}{\mathbb{Z}}

\newcommand{\cov}{\mathrm{cov}}
\newcommand{\var}{\mathrm{var}}

\newcommand{\Ex}{\mathbb{E}}

\DeclareMathOperator{\cspn}{\overline{\mathrm{sp}}}

\title{On the asymptotic behavior of a finite section of the optimal causal filter}

\author[]{Junho Yang  \thanks{\textit{Email}: \url{junhoyang@stat.sinica.edu.tw}}
}
\affil[]{Institute of Statistical Science, Academia Sinica}

\date{\today}

\begin{document}

\maketitle

\begin{abstract}
We establish an $L_1$-bound between the coefficients of the optimal causal filter applied to the data-generating process and its finite sample approximation. Here, we assume that the data-generating process is a second-order stationary time series with either short or long memory autocovariances. To derive the $L_1$-bound, we first provide an exact expression for the coefficients of the causal filter and their approximations in terms of the absolute convergent series of the multistep ahead infinite and finite predictor coefficients, respectively. Then, we prove a so-called uniform Baxter's inequality to obtain a bound for the difference between the infinite and finite multistep ahead predictor coefficients in both short and long memory time series. The $L_1$-approximation error bound for the causal filter coefficients can be used to evaluate the performance of the linear predictions of time series through the mean squared error criterion.

\vspace{0.5em}


\end{abstract}


\section{Introduction}\label{sec:intro}
Let $\{X_t\} = \{X_t\}_{t\in\Z}$ be a univariate centered second-order stationary process (the extension to the multivariate case will be discussed in Section \ref{sec:multi}). Following the terminology of \cite{p:wil-mce-16}, given filter coefficients $\{h_k\}_{k \in \Z}$, a \textit{target} random variable is defined as the filtered process of $\{X_t\}$, which can be written as
\begin{equation} \label{eq:target}
Y_n := H(B) X_n = \sum_{k \in \Z} h_k X_{n-k}, \quad n\in\N,
\end{equation}  where $B$ is the backshift operator and $H(e^{i\omega}) = \sum_{k \in \Z} h_k e^{ik\omega}$ is the frequency response function of $\{h_k\}_{k\in \Z}$. 
Two typical examples of linear filters used for either forecasting or (frequency domain) feature extraction are
$H_{m}(B)=B^{-m}$ for $m \in \N$ and $H_{\text{fs}}(e^{i\omega}) = I(\mu_1 \leq |\omega| < \mu_2)$ for some $0 \leq \mu_1 < \mu_2 \leq \pi$, respectively, where $I(\cdot)$ is an indicator function. The former shifts the original signal forward and is used in the prediction of time series (see \cite{b:pou-01} and the references therein). The latter is known as the (ideal) frequency-selective filter, which includes the well-known low-pass, high-pass, band-pass, and band-stop filters in signal processing literature (cf. \cite{b:opp01}, Section 5.1.1). For further examples of linear filters and their corresponding target random variables, we refer readers to \cite{p:wil-mce-16}, Section 2.1.

Given a target random variable $Y_n$, we are interested in finding the best linear predictor of $Y_n$ based on the infinite past and present data (often referred to as the \textit{real-time} data). Using the time-invariant property of second-order stationary processes,
the best linear predictor of $Y_n$ given $\{X_t\}_{t \leq n}$ can be expressed as
\begin{equation} \label{eq:Yhat}
\widehat{Y}_n := \sum_{k=1}^{\infty}\hat{h}_{k} X_{n+1-k}, \quad n \in \N,
\end{equation} where $\{\hat{h}_k\}_{k=1}^{\infty}$ denotes the coefficients that minimize the mean squared prediction error (MSPE)
\begin{equation*}
(\hat{h}_1, \hat{h}_2, ...) = \text{argmin}_{(a_1, a_2, ...) \in \mathbb{R}^\N}\|Y_n - \sum_{k=1}^{\infty}a_{k} X_{n+1-k}\|_V^2.
\end{equation*}
Here, $V$ is the real Hilbert space spanned by $\{X_t\}$ in $L^2(\Omega, \mathcal{F}, P)$ equipped with the norm $\Vert x\Vert_{V}:=\{\Ex[x^2]\}^{1/2}$. Techniques to solve (\ref{eq:Yhat}) for (the frequency response function of) $\{\hat{h}_k\}_{k=1}^{\infty}$ were first developed by \cite{p:hop-wie-31} using spectral decompositions (see also \cite{b:wie-49}) and more recently by \cite{p:sub-yan-23} using pseudo prediction and de-convolution techniques. \cite{p:wil-mce-16} established a general framework for prediction and model fitting problems by obtaining $\{\hat{h}_k\}_{k=1}^{\infty}$ for specific linear filters. The frequency response function of $\{\widehat{h}_{k}\}_{k=1}^{\infty}$ is also called the \textit{optimal causal (concurrent) filter} or \textit{causal Wiener filter} of $Y_n$, which has been extensively studied in time series and signal processing literature (see \cite{b:and-12}, Chapter 9 for a comprehensive review). 

In practical scenarios, only a finite number of past observations is available. Consequently, the coefficients $\{\hat{h}_k\}_{k=1}^{\infty}$ and the linear  prediction $\widehat{Y}_n$ need to be adjusted based on these finite samples (which is more realistic). Let $\widehat{Y}_{n}^{(n)}$ be the best linear predictor of $Y_n$ given the finite past observations $\{X_k\}_{k=1}^{n}$. Utilizing a similar representation as in (\ref{eq:Yhat}), $\widehat{Y}_{n}^{(n)}$ can be written as  
\begin{equation} \label{eq:Yhatn}
\widehat{Y}_{n}^{(n)} := \sum_{k=1}^{n} \hat{h}_{k,n} X_{n+1-k}, \quad n \in \N,
\end{equation}
where $\{\hat{h}_{k,n}\}_{k=1}^{n}$ are the best ``finite'' predictor coefficients of $Y_n$, which minimize the MSPE $\Ex |Y_n - \sum_{k=1}^{n} a_{k,n} X_{n+1-k}|^2$. The coefficients $\{\hat{h}_{k,n}\}_{k=1}^{n}$ are a solution to a finite-order Toeplitz system, which can be solved effectively using a Levinson-Durbin algorithm (\cite{p:lev-47} and \cite{p:dur-60}) in $O(n^{2})$ time complexity.

Then, a natural question that arises in this context is:
\begin{equation*}
\textit{How closely do the coefficients $\{\hat{h}_k\}$ approximate $\{\hat{h}_{k,n}\}$?}
\end{equation*}
 In detail, our primary objective in this article is to derive the rate of convergence to zero of the $L_1$-difference $\sum_{k=1}^{n} |\hat{h}_k - \hat{h}_{k,n}|$. Investigation of this rate is important in the prediction theory of time series literature as it is closely related to the quality of linear predictions. Let 
\begin{equation} \label{eq:Ytilde}
\widetilde{Y}_n^{(n)} = \sum_{k=1}^{n} \hat{h}_{k} X_{n+1-k}, \quad n \in \N,
\end{equation} be the $n$-th partial sum of $\widehat{Y}_n$ which belongs to the space of linear span of $X_n, ..., X_1$. Then, both $\widetilde{Y}_n^{(n)} $ and $\widehat{Y}_n $ are the approximations of $\widehat{Y}_n^{(n)}$. To assess the performance of predictions, for $n \in \N$, we consider the following (square-root of the) MSPEs:
\begin{equation} \label{eq:sigmatilde}
\widetilde{\sigma}_n := \| \widehat{Y}_n^{(n)} - \widetilde{Y}_n^{(n)}\|_V
\end{equation} and
\begin{equation}\label{eq:sigma}
\sigma_n := \| \widehat{Y}_n^{(n)} - \widehat{Y}_n\|_V.
\end{equation}
Using a triangular inequality, it is easily seen that for all $n\in \N$,
\begin{eqnarray}
\widetilde{\sigma}_n  =   \bigg\| \sum_{k=1}^{n} \{ \hat{h}_k - \hat{h}_{k,n}\} X_{n+1-k} \bigg\|_V
&\leq& \sum_{k=1}^{n}  | \hat{h}_k - \hat{h}_{k,n}| \|X_{n+1-k}\|_V \nonumber \\
&=& \sqrt{\var X_0} \sum_{k=1}^{n} |\hat{h}_k - \hat{h}_{k,n}| 
\label{eq:sigmatilde-bound}
\end{eqnarray} and
\begin{eqnarray}
\sigma_n &=& \bigg\| \sum_{k=1}^{n} \{ \hat{h}_k - \hat{h}_{k,n}\} X_{n+1-k} +  \sum_{k=n+1}^{\infty} \hat{h}_k X_{n+1-k} \bigg\|_V \nonumber \\
&\leq& \sqrt{\var X_0} \sum_{k=1}^{n} |\hat{h}_k - \hat{h}_{k,n}| 
+ \bigg\| \sum_{k=n+1}^{\infty} \hat{h}_k X_{n+1-k} \bigg\|_V.
\label{eq:sigma-bound}
\end{eqnarray} 
Therefore, the convergence rate of $\sum_{k=1}^{n} |\hat{h}_k - \hat{h}_{k,n}|$ as $n\rightarrow \infty$ plays an important role in the asymptotic behavior of $\widetilde{\sigma}_n $ and $\sigma_n$. 

\subsection{Literature review and contributions}
In this article, we derive an exact asymptotic rate for $\sum_{k=1}^{n} |\hat{h}_k - \hat{h}_{k,n}|$, $\widetilde{\sigma}_n$, and $\sigma_n$ under both short and long memory frameworks (refer to Assumptions \ref{assum:SM} and \ref{assum:LM}, respectively). These results are a good improvement in the existing literature, which we will discuss in greater detail below.

\vspace{0.5em}

\noindent \textbf{Short memory framework:} 
Under the short memory framework, that is $\sum_{r\in \Z} |\cov(X_r,X_0)|<\infty$,
\cite{p:bax-62} derived an asymptotic rate for $\sum_{k=1}^{n} |\hat{h}_k - \hat{h}_{k,n}|$ and $\widetilde{\sigma}_n$ when $H(B)=B^{-1}$ (see also \cite{p:dev-64, p:ibr-65} and the survey paper by \cite{p:bin-12}). In \cite{B63}, he proved a norm inequality for a finite-section Wiener--Hopf equation, which can be stated in our context as $\sum_{k=1}^{n} |\hat{h}_k - \hat{h}_{k,n}| = O(\sum_{k=n+1}^{\infty} |\hat{h}_k|)$ as $n\rightarrow \infty$. Here, we require a condition that $\{\hat{h}_k\}$ is absolutely summable. However, this summability condition is unwieldy and the tail convergence rate of $\{\hat{h}_k\}$ is not transparent.

In this article, we extend the existing literature on short memory processes in at least two ways. Firstly, provided that $\{h_k\}$ is absolutely summable, we derive an exact asymptotic rate for 
 $\sum_{k=1}^{n} |\hat{h}_k - \hat{h}_{k,n}|$, $\widetilde{\sigma}_n$, and $\sigma_n$ in terms of the tail behavior of $\{h_k\}$ and the decay rate of the infinite-order MA and AR coefficients of $\{X_t\}$ (see Assumption \ref{assum:H}(H-SM1)). The results are summarized in Theorems \ref{thm:L1-short}(i) and \ref{thm:sigma-SM}(i). We mention that we prove these results through the lens of linear prediction, which is distinct from the orthogonal polynomial technique that is used in majority of the existing works, including \cite{p:bax-62,B63}. Secondly, under stronger summability conditions such as $\sum_{r\in \Z} (1+|r|) |\cov(X_r,X_0)| <\infty$ (see Assumption \ref{assum:H}(H-SM2) and Remark \ref{rmk:cond}(i)), our convergence results can be extended to square-summable $\{h_k\}$ (Theorems \ref{thm:L1-short}(ii) and \ref{thm:sigma-SM}(ii)). As a result, we are able to derive the convergence results for frequency-selective filters which have square summable filter coefficients (but not absolutely summable). 

Comparing our asymptotic rates with existing results, when $H(B) = B^{-1}$ (or, more generally, when $H(\cdot)$ corresponds to a finite impulse response (FIR) filter), our asymptotic rate for $\sum_{k=1}^{n} |\hat{h}_k - \hat{h}_{k,n}|$ and $\widetilde{\sigma}_n$ coincides with those in \cite{p:bax-62} (see Remarks \ref{rmk:fir1} and \ref{rmk:fir2}). However, when $H(\cdot)$ corresponds to a general linear filter, an additional tail sum of $\{h_k\}$ is introduced in our asymptotic error bounds.

\vspace{0.5em}

\noindent \textbf{Long memory framework:}
Due to the technical challenges, there is very little work on the convergence rates of $\sum_{k=1}^{n} |\hat{h}_k - \hat{h}_{k,n}|$, $\widetilde{\sigma}_n$, and $\sigma_n$ so far under the long memory framework, that is, $\sum_{r\in \Z} |\cov(X_r, X_0)| =\infty$. A remarkable exception in this area is \cite{p:ino-kas-06}, where they investigated the asymptotic rate of $\sum_{k=1}^{n} |\hat{h}_k - \hat{h}_{k,n}|$ for $H(B) = B^{-1}$ under the long memory framework.

In this article, we assume that the infinite-order MA and AR coefficients of $\{X_t\}$ satisfy condition (LM$_d$) for some $d \in (0,1/2)$, where $d$ is a fractional differencing parameter of $\{X_t\}$ (see Assumption \ref{assum:LM} and Remark \ref{rmk:cond}(iv)). Therefore, $\{X_t\}$ exhibits long memory phenomenon (Remark \ref{rmk:cond}(iii)). Under this condition together with filter coefficients satisfying $\sum_{k \in \Z} (1+|k|)^d |h_k| <\infty$, we derive the convergence rate of $\sum_{k=1}^{n} |\hat{h}_k - \hat{h}_{k,n}|$, $\widetilde{\sigma}_n$, and $\sigma_n$ (Theorems \ref{thm:L1-long} and \ref{thm:sigma-LM}). We mention that when $H(B)=B^{-1}$, our asymptotic rate for $\sum_{k=1}^{n} |\hat{h}_k - \hat{h}_{k,n}|$ aligns with that in \cite{p:ino-kas-06}. However, our results extend to more general linear filters, including filter coefficients of (i) finite impulse response, (ii) geometric decay, or (iii) polynomial decay, where $|h_k| \leq (1+|k|)^{-(1+d+\epsilon)}$ for some $\epsilon>0$.

\vspace{0.5em}

\noindent \textbf{Uniform Baxter's inequality:}
A crucial part in proving our main asymptotic results is on the stronger convergence results of the multi-step ahead predictor coefficients. Below, we explain this in more detail.

For $m,n \in \N$, let $\{\phi_k^m\}_{k=1}^{\infty}$ and $\{\phi_{k,n}^{m}\}_{k=1}^{n}$ denote the solution of (\ref{eq:Yhat}) and (\ref{eq:Yhatn}), respectively, corresponding to the frequency response function of form $H_m(B) = B^{-m}$. The coefficients $\{\phi_k^m\}_{k=1}^{\infty}$ (resp., $\{\phi_{k,n}^{m}\}_{k=1}^{n}$) are referred to as the $m$-step ahead infinite (resp., finite) predictor coefficients of $\{X_t\}$ (see Section \ref{sec:prelim} for more details). As mentioned previously, existing works such as \cite{p:bax-62} (under short memory) and \cite{p:ino-kas-06} (under long memory) proved convergence results for $\sum_{k=1}^{n} |\phi_k^1 - \phi_{k,n}^1|$ with rate of convergence to zero, while \cite{p:fin-91} extended the results to arbitrary fixed $m \in \N$ under short memory framework. Therefore, these works focus on the ``point-wise'' convergence results in the sense that $m$ is fixed and $n$ approaches infinity.

In this article, we prove a uniform-type convergence of $\sum_{k=1}^{n} |\phi_k^m - \phi_{k,n}^m|$ under both short and long memory frameworks, which, as far as we know, is the first attempt. The precise statements of these results can be found in Theorems \ref{thm:unif-BaxterSM} and \ref{thm:unif-BaxterLM}. We refer to these results as a \textit{uniform Baxter's inequality} (see also the prototype inequality in (\ref{eq:Baxter})). The uniform Baxter's inequality for both short and long memory processes is not only a special case of our results but also a key ingredient in proving our main asymptotic results. Furthermore, the two versions of uniform Baxter's inequality may also be of independent interest in time series literature.

\subsection{Organization}
The rest of the article is organized as follows. In Section \ref{sec:prelim}, we introduce most of the notation that is used throughout the article and present basic results. In Section \ref{sec:main-results}, we state the main results. In Section \ref{sec:unif-Baxter}, we state and prove a uniform Baxter's inequality for both short and long memory time series. In Section \ref{sec:proof-main}, we prove our main results. In Section \ref{sec:multi}, we discuss the multivariate generalization. Finally, technical lemmas can be found in Appendix \ref{appen:A}.


\section{Preliminaries} \label{sec:prelim}

Let $\{X_t\} = \{X_t\}_{t\in\Z}$ be a univariate centered second-order stationary process defined on the probability space $(\Omega, \mathcal{F}, P)$.
The corresponding autocovariance and spectral density function of $\{X_t\}$ are denoted by $\gamma(k) = \cov (X_k, X_0)$ ($k \in \mathbb{Z}$) and
$f(e^{i\omega}) = \sum_{k \in \mathbb{Z}} \gamma(k) e^{ik\omega}$ ($\omega \in [0,2\pi)$), respectively. 
Let $V$ be the real Hilbert space spanned by $\{X_t\}$ in $L^2(\Omega, \mathcal{F}, P)$ with the inner product $(x, y)_{V}:=\Ex[xy]$ and
norm $\Vert x\Vert_{V}:=(x,x)_{V}^{1/2}$ for $x,y \in V$. For $K\subset \Z$, such as 
$(-\infty,n]:=\{n,n-1,...\}$, $[n,\infty):=\{n,n+1,...\}$,
and $[m,n]:=\{m,...,n\}$ with $m\le n$,
we define the closed subspace $V_K^{}$ of $V$ as
$V_K^{}:=\cspn \{X_k:  k\in K\}$ and $P_K$ denotes the orthogonal projection operators of $V$ onto $V_K^{}$.
Therefore, $\widehat{Y}_n$ and $\widehat{Y}_{n}^{(n)}$ in (\ref{eq:Yhat}) and (\ref{eq:Yhatn}) can be expressed as
\begin{equation*}
\widehat{Y}_n = P_{(-\infty,n]}Y_n \quad \text{and} \quad 
\widehat{Y}_{n}^{(n)} = P_{[1,n]}Y_n, \quad n \in \N.
\end{equation*} 
For $p \in [1,\infty)$, $\ell_p(K)$ denotes the space of all real-valued sequences $(a_k: k \in K)$ such that $\sum_{k\in K} |a_k|^p<\infty$. 
We write $\ell_p :=\ell_p(\Z)$, $\ell_p^{+} := \ell_p(\N\cup \{0\})$, and $\ell_p^{-} = \ell_p(\{-1,-2, ....\})$. Lastly, 
for sequences $\{a_j\}_{j=0}^{\infty}$ and $\{b_j\}_{j=0}^{\infty}$, $a_j \sim b_j$ indicates that $\lim_{j\rightarrow \infty} a_j/b_j = 1$.

To set down Wiener's prediction theory, throughout this paper we assume that $\{X_t\}$ is purely nondeterministic. That is,
$\bigcap_{n = -\infty}^{\infty} V_{(-\infty,n]} = \{0\}$. 
It is well-known that $\{X_t\}$ is purely nondeterministic if and only if  $f(e^{i\omega})$ exists and $\int_{0}^{2\pi} |\log f(e^{i\omega})| d\omega < \infty$ (cf. \cite{b:roz-67}, Section 4). Since $\{X_t\}$ is purely nondeterministic, we apply the Szeg\"{o}'s identity (\cite{p:sze-21, p:wie-mas-58})
\begin{equation} \label{eq:spec-decom1}
f(e^{i\omega}) = \sigma^2|\psi(e^{i\omega})|^2 = \sigma^2 \big| 1+ \sum_{j=1}^{\infty} \psi_{j} e^{-ij\omega} \big|^2, \quad \omega \in [0,2\pi).
\end{equation} Here, $\{\psi_j\}_{j=0}^{\infty}$ (setting $\psi_0=1$) is a real-valued sequence in $\ell_2^+$. We note that the characteristic polynomial
$\Psi(z) = 1+ \sum_{j=1}^{\infty} \psi_{j} z^j$ has no roots inside the unit circle. 
The coefficients $\{\psi_j\}$ are known as the infinite-order moving average (MA) coefficients of $\{X_t\}$ since 
these coefficients appear in the classical Wold's decomposition (\cite{b:doo-53}) of $\{X_t\}$. Specifically, $X_t$ can be represented as
\begin{eqnarray} \label{eq:MAinfty}
X_t = \sum_{j=0}^{\infty} \psi_{j} \varepsilon_{t-j}, \quad t \in \mathbb{Z},
\end{eqnarray} where $\{\varepsilon_t\}$ is a white noise process with variance $\sigma^2>0$ that satisfies 
$\cspn \{\varepsilon_k:  k \leq t\} = \cspn \{X_k:  k \leq t\}$ for all $t \in \mathbb{Z}$.

Since $\Psi(z)$ in (\ref{eq:spec-decom1}) has no roots inside the unit circle, we can define the autoregressive (AR) transfer function as 
\begin{equation} \label{eq:ARfunction}
\Phi(z) = \Psi(z)^{-1} = 1 - \sum_{j=1}^{\infty} \phi_j z^j, \quad |z| < 1.
\end{equation} Here, $\{\phi_j\}_{j=0}^{\infty}$ (setting $\phi_0=-1$) is also a real-valued sequence in $\ell_2^+$.
Under the additional assumption that $\{\phi_j\} \in \ell_1^+$,
$\{X_t\}$ admits the Wold-type one-sided infinite-order AR representation
\begin{equation}\label{eq:ARinfty}
X_t - \sum_{j=1}^{\infty} \phi_{j}X_{t-j} =  \varepsilon_{t},  \quad t \in \mathbb{Z},
\end{equation} where $\{\varepsilon_t\}$ is the same white noise process in (\ref{eq:MAinfty}) (see \cite{b:pou-01}, Sections 5-6, and \cite{p:kre-11}, page 706). Due to the above representation, we refer to $\{\phi_j\}$ as the infinite-order AR coefficients of $\{X_t\}$.

Finally, we introduce the multistep ahead infinite and finite predictor coefficients. These coefficients appear in the expression of $\{\hat{h}_k\}$ and $\{\hat{h}_{k,n}\}$ (see Theorem \ref{thm:expression}) and are the main objects in the uniform Baxter's inequality in Section \ref{sec:unif-Baxter}.
 For a linear filter of form $H_m(B) = B^{-m}$, $m \in \N$, the corresponding target random variable is $Y_n  = H_m(B)X_n = X_{n+m}$, $n \in \N$. Then, the best linear predictor of $Y_n = X_{n+m}$ given $\{X_{t}\}_{t\leq n}$ is
\begin{equation} \label{eq:m-step1}
\widehat{Y}_n = \widehat{X}_{n+m} = P_{(-\infty,n]} X_{n+m} = \sum_{k=1}^{\infty} \phi^{m}_k X_{n+1-k},
\end{equation} where $\{ \phi^{m}_k\}_{k=1}^{\infty}$ are the $m$-step ahead infinite predictor coefficients (see, e.g., \cite{b:bro-dav-06}, Chapter 5). For example, by using \cite{SRY}, Lemma A.2, $\phi^{m}_k$ and the infinite-order MA and AR coefficients of $\{X_t\}$ are related by the following identity:
\begin{equation} \label{eq:BLPinfty_ARMA1}
\phi_{k}^m = \sum_{\ell=0}^{m-1} \psi_{\ell} \phi_{k+m-1-\ell}, \quad m,k \in \N.
\end{equation} A linear prediction of $Y_n$ based on the finite observations $X_n, ..., X_1$ is
\begin{equation} \label{eq:m-step2}
\widehat{Y}_n^{(n)} = \widehat{X}_{n+m}^{(n)} = P_{[1,n]} X_{n+m} = \sum_{k=1}^{n} \phi^{m}_{k,n} X_{n+1-k}.
\end{equation} The coefficients $\{\phi^{m}_{k,n}\}_{k=1}^{n}$ are referred to as the $m$-step ahead finite predictor coefficients.

\section{Main results} \label{sec:main-results}

In this section, we obtain bounds for $\sum_{k=1}^{n} |\hat{h}_k - \hat{h}_{k,n}|$, $\widetilde{\sigma}_n$, and $\sigma_n$.
 To do so, we first define the space of sequences with fast-decaying coefficients. For $\alpha \in [0,\infty)$ and $K \subset \Z$, let
\begin{equation} \label{eq:Calpha}
\mathcal{C}_\alpha (K) = \{ (a_k: k\in K):~ \sum_{k\in K} (1+|k|)^\alpha |a_k| < \infty \}.
\end{equation} We write $\mathcal{C}_\alpha := \mathcal{C}_{\alpha}(\Z)$, $\mathcal{C}_{\alpha}^{+} := \mathcal{C}_{\alpha}(\N\cup \{0\})$, and $\mathcal{C}_{\alpha}^{-} := \mathcal{C}_{\alpha}(\{-1,-2, ....\})$. 
Next, a positive measurable function $\ell(x)$ defined on $[0,\infty)$ is called a \textit{slowly varying} function, if $\lim_{x\rightarrow \infty} \ell(\lambda x) / \ell(x) = 1$ for all $\lambda >0$ (cf. \cite{b:bin-89}, Chapter 1).

With this notation, we consider the following two scenarios for the infinite-order MA and AR coefficients of $\{X_t\}$.

\begin{assumption}[Short memory time series] \label{assum:SM}
$\{X_t\}$ is purely nondeterministic and the infinite-order MA and AR coefficients of $\{X_t\}$ are such that
\begin{equation} \tag{SM$_\alpha$}
\{\psi_j\}, \{\phi_j\} \in \mathcal{C}_\alpha^{+} \quad \text{for some}~~ \alpha \in [0,\infty).
\end{equation}
\end{assumption}

\begin{assumption}[Long memory time series] \label{assum:LM}
$\{X_t\}$ is purely nondeterministic and for some $d\in (0,1/2)$ and slowly varying function $\ell(\cdot)$,
the infinite-order MA and AR coefficients of $\{X_t\}$ are such that
\begin{equation} \tag{LM$_d$}
\psi_j \sim j^{-(1-d)} \ell(j) \quad \text{and} \quad \phi_j \sim j^{-(1+d)} \frac{1}{\ell(j)} \frac{d\sin(\pi d)}{\pi}
\quad \text{as}~~ j \rightarrow \infty.
\end{equation}
\end{assumption}
We make some remarks on the conditions above.
\begin{remark} \label{rmk:cond}
\begin{itemize}
\item[(i)] Under (SM$_\alpha$), the spectral density and autocovariance functions are such that $0<\inf_\omega f(e^{i\omega}) \leq \sup_\omega f(e^{i\omega}) <\infty$ and $\sum_{k\in\Z} |\gamma(k)| <\infty$. Thus, $\{X_t\}$ has a short memory. A sufficient condition for (SM$_\alpha$) to hold is that $\inf_{\omega} f(e^{i\omega}) >0$, and the autocovariance function $\{\gamma(k)\}$ belongs to $\mathcal{C}_\alpha$ (see, e.g., \cite{p:kre-11}, Lemma 2.1). 

\item[(ii)] For $\alpha^{\prime} \geq 0$, let
$\|\psi\|_{\alpha^\prime} = \sum_{j=0}^{\infty} (1 + j)^{\alpha^\prime} |\psi_j|$ and
$\|\phi\|_{\alpha^\prime} = \sum_{j=0}^{\infty} (1 + j)^{\alpha^\prime} |\phi_j|$. Then, under (SM$_\alpha$) 
for some $\alpha \geq 0$ and for $\alpha^\prime \in [0,\alpha]$, $\|\psi\|_{\alpha^\prime}$ and $\|\psi\|_{\alpha^\prime}$ are both finite.
These quantities appear in the constant factor of the uniform Baxter's inequality for short memory processes (see Theorem \ref{thm:unif-BaxterSM} below).

\item[(iii)] Using \cite{p:ino-kas-06}, equation (2.22), under (LM$_d$), we have
\\ $\gamma(k) \sim k^{-(1-2d)}\ell(k)^2 B(d,1-2d)$ as $k \rightarrow \infty$, where $B(\cdot, \cdot)$ is the beta function. This indicates that for $d\in (0,1/2)$, $\sum_{k\in\Z} |\gamma(k)|=\infty$. Thus, $\{X_t\}$ has a long memory. 

\item[(iv)] An autoregressive fractional integrated moving average (ARFIMA) process (\cite{p:gra-joy-80, p:hos-81}) with the fractional difference parameter $d \in (0,1/2)$ satisfies the condition (LM$_d$) for some constant function $\ell(\cdot)$
(see \cite{p:kok-95}, equation (2.2) and Corollary 3.1).
\end{itemize}
\end{remark}
Lastly, to obtain the well-defined target random variable, we assume the following for the filter coefficients.
\begin{assumption} \label{assum:H}
A time series $\{X_t\}$ and the filter coefficients $\{h_k\}$ satisfy one of the following three conditions:
\begin{equation} \tag{H-SM1}
\text{$\{X_t\}$ satisfies (SM$_\alpha$) for some $\alpha \geq 0$ and } \{h_k\} \in \ell_1.
\end{equation}
\begin{equation} \tag{H-SM2}
\text{$\{X_t\}$ satisfies (SM$_\alpha$) for some $\alpha \geq 1$ and } \{h_k\} \in \ell_2.
\end{equation}
\begin{equation} \tag{H-LM}
\text{$\{X_t\}$ satisfies (LM$_d$) for some $d \in (0,1/2)$ and } \{h_k\} \in \mathcal{C}_d.
\end{equation}
\end{assumption}
Under one of the three conditions that are stated above, it is easily seen that the target random variable $Y_n$ in (\ref{eq:target}) is well-defined on $V$.


In the theorem below, we show that $\{\hat{h}_k\}_{k=1}^{\infty}$ and $\{\hat{h}_{k,n}\}_{k=1}^{n}$ can be expressed in terms of the absolute convergent series of  the multistep ahead infinite and finite predictor coefficients, respectively.
\begin{theorem} \label{thm:expression}
Suppose $\{X_t\}$ satisfies either Assumption \ref{assum:SM} or \ref{assum:LM} and the filter coefficients $\{h_k\}$ satisfy Assumption \ref{assum:H}. Let $\{\phi_{k}^{m}\}_{k=1}^{\infty}$ and $\{\phi_{k,n}^{m}\}_{k=1}^{n}$ are the $m$-step ahead infinite and finite predictor coefficients, respectively (see (\ref{eq:m-step1}) and (\ref{eq:m-step2})). 
Then,
\begin{equation} \label{eq:FIRcoeff1}
\hat{h}_k = \sum_{m=1}^{\infty}h_{-m}  \phi_{k}^m  + h_{k-1}, \quad k \in \N,
\end{equation} 
and
\begin{equation} \label{eq:FIRcoeff2}
\hat{h}_{k,n} = \sum_{m=1}^{\infty} h_{-m} \phi_{k,n}^m  +h_{k-1} +  \sum_{m=1}^{\infty} h_{n-1+m} \phi_{n+1-k,n}^m
, \quad n \in \N \quad \text{and} \quad k\in\{1,...,n\}.
\end{equation} Here, $\{h_{-k}\}_{m=1}^{\infty}$ corresponds to the filter coefficients of the ``future'' samples $X_{n+1}, X_{n+2}, \cdots$.
For sufficiently large enough $n \in \N$, the above sums converge absolutely. 
Furthermore, under (H-SM1) or (H-LM),
\begin{equation} \label{eq:hatH1}
\{\hat{h}_k\}_{k=1}^{\infty} \in \ell_1(\N)
\end{equation}  and under (H-SM2),
\begin{equation}  \label{eq:hatH2}
\{\hat{h}_k\}_{k=1}^{\infty} \in \ell_2(\N).
\end{equation}
\end{theorem}
The corollary below immediately follows from the above theorem.

\begin{corollary} \label{coro:expression}
Suppose the same set of assumptions and notation as in Theorem \ref{thm:expression} holds. Then, for large enough $n\in \N$,
\begin{equation}
 \sum_{k=1}^{n} | \hat{h}_k - \hat{h}_{k,n}| 
\leq  \sum_{m=1}^{\infty}  | h_{-m}|  \sum_{k=1}^{n} | \phi_{k}^m - \phi_{k,n}^m|
+  \sum_{m=1}^{\infty} \sum_{k=1}^{n}  |h_{n+m-1}|   |\phi_{n+1-k,n}^{m}|.
 \label{eq:hatHdiff}
\end{equation}
\end{corollary}

From the above corollary, the $L_1$-difference $\sum_{k=1}^{n} | \hat{h}_k - \hat{h}_{k,n}|$ can be bounded using two terms. The bound for the first term requires a uniform Baxter's inequality, which will be discussed in detail in Section \ref{sec:unif-Baxter}. The error bound for the second term in (\ref{eq:hatHdiff}) is provided in Appendix, Lemma \ref{lemma:abs-conv2}. 
By combining these two results, we state the $L_1$-approximation theorem for short memory time series.

\begin{theorem}[$L_1$-approximation theorem for short memory time series] \label{thm:L1-short}
Suppose $\{X_t\}$ is a short memory time series that satisfies Assumption \ref{assum:SM} for some $\alpha \geq 0$.
Then, the following two assertions hold.
\begin{itemize}
\item[(i)] If the filter coefficients $\{h_k\}$ belong to $\ell_1$, then
\begin{equation} \label{eq:short-baxter11}
\sum_{k=1}^{n} |\hat{h}_k - \hat{h}_{k,n}| = o(n^{-\alpha}) + O\big(\sum_{k=n}^{\infty} |h_k|\big)
 \quad \text{as}~~ n \rightarrow \infty.
\end{equation}
\item[(ii)]  If $\alpha \geq 1$ and the filter coefficients $\{h_k\}$ belong to $\ell_2$, then
\begin{equation} \label{eq:short-baxter22}
\sum_{k=1}^{n} |\hat{h}_k - \hat{h}_{k,n}| = o(n^{-\alpha+1}) + O\big(\sup_{k \geq n} |h_k|\big)
\quad \text{as}~~ n \rightarrow \infty.
\end{equation}
\end{itemize}
\end{theorem}

\begin{remark} \label{rmk:constant-factor}
As pointed out by the Associate Editor, in certain applications, not only the rates of the convergence to zero but also constant factors are of interest. Precise expressions of the constant factors, which are multiplied by the big-$O$ and little-$o$ asymptotic notation in (\ref{eq:short-baxter11}) and (\ref{eq:short-baxter22}), can be found in (\ref{eq:SM1-exactbound}) and (\ref{eq:SM2-exactbound}), respectively. 
\end{remark}

\begin{remark} \label{rmk:fir1}
We provide a brief discussion of the above bounds for some special cases.
\begin{itemize}
\item[(i)] Suppose $H(\cdot)$ corresponds to an FIR filter. In this case, $\{h_k\} \in \ell_1$ and the term $O\left(\sum_{k=n}^{\infty} |h_k|\right)$ in (\ref{eq:short-baxter11}) vanishes for some large enough $n \in \N$. Therefore, we have $\sum_{k=1}^{n} |\hat{h}_k - \hat{h}_{k,n}| = o(n^{-\alpha})$ as $n \rightarrow \infty$. This rate coincides with that in \cite{p:bax-62}, Theorem 2.2.

\item[(ii)] Suppose the MA and AR coefficients of $\{X_t\}$ decay geometrically fast to zero (e.g., ARMA processes) and the filter coefficients $\{h_k\}$ belong to $\ell_2$. Then, there exists $\rho \in (0,1)$ such that $\sum_{k=1}^{n} |\hat{h}_k - \hat{h}_{k,n}| = O(\rho^n) + O(\sup_{k \geq n} |h_k|)$ as $n \rightarrow \infty$. A proof of this bound is similar to that in the proof of (\ref{eq:short-baxter22}), thus we omit the details.
\end{itemize}
\end{remark}

Below is the analogous $L_1$-approximation theorem for long memory time series.

\begin{theorem}[$L_1$-approximation theorem for long memory time series] \label{thm:L1-long}
Suppose $\{X_t\}$ is a long memory time series that satisfies Assumption \ref{assum:LM} for some $d \in (0,1/2)$. Furthermore, we assume that the filter coefficients $\{h_k\}$ belong to $\mathcal{C}_d$. Then,
\begin{equation} \label{eq:long-baxter11}
\sum_{k=1}^{n} |\hat{h}_k - \hat{h}_{k,n}| = O(n^{-d}) \quad \text{as}~~ n \rightarrow \infty.
\end{equation}
\end{theorem}
We mention that in the special case where $H(B) = B^{-1}$, the bound in (\ref{eq:long-baxter11}) coincides with that in \cite{p:ino-kas-06}, Theorem 4.1. The precise constant factor that is multiplied by the big-$O$ asymptotic notation in (\ref{eq:long-baxter11}) can be found in (\ref{eq:LM-exactbound}).

\begin{remark}[Estimation of the finite predictor coefficients]
Even though the coefficients $\{\hat{h}_{k,n}\}$ are constructed in the spirit of finite sample observations, they are functions of the unknown autocovariances of $\{X_t\}$. Therefore, $\{\hat{h}_{k,n}\}$ cannot be directly evaluated. To estimate the low-dimensional coefficients $\{\hat{h}_{k,p}\}_{k=1}^{p}$ based on the ``full'' observations $X_n, ..., X_1$, where $p = p(n)$ is such that $1/p + p/n \rightarrow 0$ as $n \rightarrow \infty$, we first estimate the multistep ahead finite predictor coefficients by fitting AR$(p)$ model to the data (the so-called plug-in method; see \cite{p:bha-96}). Then, we estimate the coefficients $\{\hat{h}_{k,p}\}_{k=1}^{p}$ by truncating equation (\ref{eq:FIRcoeff2}) and replacing the finite predictor coefficients with their sample estimators. Probabilistic bounds that are developed in \cite{p:kre-11,SRY} (under short memory) and \cite{p:pos-07} (under long memory) can be used to obtain the sampling properties of the resulting estimator. However, we leave this investigation for future research.
\end{remark}

Next, we consider the two (square-root of the) MSPEs: $\widetilde{\sigma}_n$ and $\sigma_n$ that are defined in (\ref{eq:sigmatilde}) and (\ref{eq:sigma}), respectively. The following theorem provides bounds for $\widetilde{\sigma}_n$ and $\sigma_n$ in the case of short memory time series.

\begin{theorem}[MSPEs for short memory time series]
\label{thm:sigma-SM}
Suppose $\{X_t\}$ is a short memory time series that satisfies Assumption \ref{assum:SM} for some $\alpha \geq 0$.
Then, the following two assertions hold.
\begin{itemize}
\item[(i)] If the filter coefficients $\{h_k\}$ belong to $\ell_1$, then
\begin{equation} \label{eq:sigma-SM1}
\widetilde{\sigma}_n = o(n^{-\alpha}) + \big(\sum_{k=n}^{\infty} |h_k|\big)
\quad \text{and} \quad \sigma_n = o(n^{-\alpha})  + \big(\sum_{k=n}^{\infty} |h_k|\big)
\quad \text{as}~~ n\rightarrow \infty.
\end{equation}

\item[(ii)]  If $\alpha \geq 1$ and the filter coefficients $\{h_k\}$ belong to $\ell_2$, then
\begin{equation} \label{eq:sigma-SM2}
\widetilde{\sigma}_n = o(n^{-\alpha+1}) + O\big(\sup_{k \geq n} |h_k|\big)
\quad \text{and} \quad \sigma_n = o(n^{-\alpha+1})  + O\big( \big\{\sum_{k=n}^{\infty} |h_k|^2\big\}^{1/2} \big)
\quad \text{as}~~ n\rightarrow \infty.
\end{equation}
\end{itemize} 
\end{theorem}

\begin{remark} \label{rmk:fir2}
Using similar arguments as in Remark \ref{rmk:fir1}, the following two assertions hold.
\begin{itemize}
\item[(i)] Suppose $H(\cdot)$ corresponds to an FIR filter. Then, both $\widetilde{\sigma}_n$ and $\sigma_n$ are $o(n^{-\alpha})$ as $n \rightarrow \infty$. These convergence rates coincide with those in \cite{p:bax-62}, Theorem 3.1.

\item[(ii)] Suppose the MA and AR coefficients of $\{X_t\}$ decay geometrically fast to zero and the filter coefficients $\{h_k\}$ belong to $\ell_2$. Then, there exists $\rho \in (0,1)$ such that 
$\widetilde{\sigma}_n = O(\rho^n) + O\big(\sup_{k \geq n} |h_k|\big)$ and $\sigma_n = O(\rho^n)  + O\big( \big\{\sum_{k=n}^{\infty} |h_k|^2\big\}^{1/2} \big)$ as $n \rightarrow \infty$. It is worth mentioning that these geometric decaying error rates were investigated in \cite{p:ibr-65} and \cite{p:pou-89} specifically for the special case when $H(B)=B^{-1}$.
\end{itemize}
 
\end{remark}

Below is the analogous results for long memory times series.
\begin{theorem}[MSPEs for long memory time series] \label{thm:sigma-LM}
Suppose $\{X_t\}$ is a long memory time series that satisfies Assumption \ref{assum:LM} for some $d \in (0,1/2)$.
Furthermore, we assume that  the filter coefficients $\{h_k\}$ belong to $\mathcal{C}_d$. Then,
\begin{equation} \label{eq:sigma-LM}
\widetilde{\sigma}_n = O(n^{-d}) \quad \text{and} \quad \sigma_n = O(n^{-d})
\quad \text{as}~~ n\rightarrow \infty.
\end{equation}
\end{theorem}

\begin{remark}
Suppose $\{X_t\}$ is an ARFIMA$(p,d,q)$ process with $d\in(0,1/2)$. Thus, $\{X_t\}$ satisfies the condition (LM$_d$). Then,
\cite{I02}, Theorem 4.3, proved that 
\begin{equation} \label{eq:sigma-new-bound}
\rho_n := \left| \|  X_{n+1} - P_{(-\infty,n]}X_{n+1}\|_V - \|  X_{n+1} - P_{[1,n]}X_{n+1}\|_V\right|  = O(n^{-1}) \quad \text{as}~~
n\rightarrow \infty.
\end{equation} Using a triangular inequality, when $H(B)=B^{-1}$, one can easily seen that $\sigma_n \geq \rho_n$, $n \in \N$. Therefore, the above is not contradictory to (\ref{eq:sigma-LM}). Nevertheless, the $O(n^{-1})$ bound in (\ref{eq:sigma-new-bound}) may suggests that we can improve the asymptotic bound for $\sigma_n$ in (\ref{eq:sigma-LM}). The proof of (\ref{eq:sigma-new-bound}) is based on a novel series expansion of $\rho_n$ in \cite{I02}, Theorem 4.1. (see also \cite{IKP2}, Theorem 5.3). On the other hand, it seems that obtaining a simple series expansion of $\sigma_n$ (or $\widetilde{\sigma}_n$) is not feasible, and determining the optimal bound for $\sigma_n$ (or $\widetilde{\sigma}_n$) remains an open question.
\end{remark}

\section{Uniform Baxter's inequality for multistep ahead finite predictor coefficients} \label{sec:unif-Baxter}

\subsection{Uniform Baxter's inequality} \label{sec:Baxter-main}

From Corollary \ref{coro:expression}, the $L_1$-difference $\sum_{k=1}^{n} | \hat{h}_k - \hat{h}_{k,n}|$ is bounded by $\sum_{m=1}^{\infty}  | h_{-m}|  \sum_{k=1}^{n} | \phi_{k}^m - \phi_{k,n}^m|$ plus an additional error due to the tail term. In this section, our objective is to obtain the bound for $\sum_{k=1}^{n}|\phi_k^m - \phi_{k,n}^m|$ for all $m \in \N$. 

According to \cite{p:bax-62}, under (SM$_\alpha$) for some $\alpha \geq 0$, there exists a constant $C \in (0,\infty)$ such that
\begin{equation} \label{eq:Baxter}
\sum_{k=1}^{n} |\phi_k^1 - \phi_{k,n}^1| \leq C \sum_{k=n+1}^{\infty} |\phi_k^1|, \quad n\in \N.
\end{equation} Therefore, $\sum_{k=1}^{n} |\phi_k^1 - \phi_{k,n}^1|$ is upper bounded by the tail sum of $\{\phi_{k}^{1}\}$, which we can easily obtain the bound in terms of the infinite-order MA and AR coefficients using (\ref{eq:BLPinfty_ARMA1}). Inequality (\ref{eq:Baxter}) is called \textit{Baxter's inequality} and by using \cite{p:fin-91}, one can extend Baxter's inequality for the $m$-step ahead finite prediction coefficients (under (SM$_\alpha$) for some $\alpha \geq 0$). Additionally, \cite{p:ino-kas-06} proved the analogous Baxter's inequality (for one-step ahead finite predictor coefficients) under (LM$_d$) for some $d \in (0,1/2)$.  

However, since the first term of the upper bound of $\sum_{k=1}^{n} | \hat{h}_k - \hat{h}_{k,n}|$ is in terms of the weighted sum of ``all'' $ \sum_{k=1}^{n} | \phi_{k}^m - \phi_{k,n}^m|$ for $m \in \N$, we require a stronger version of Baxter's inequality. In the next two theorems, we extend the results of \cite{p:bax-62} and \cite{p:ino-kas-06} and provide a \textit{uniform Baxter's inequality} under both short and long memory frameworks. 

Before stating our results, we make the constants involved in the uniform Baxter's inequality explicit. Under the short memory framework, for $\varepsilon \in (0,1)$, let
\begin{equation} \label{eq:C1}
C_1 = \frac{3-\varepsilon}{1-\varepsilon} \big\{ \sum_{j=0}^{\infty} |\psi_j| \big\}^2  \big\{ \sum_{j=0}^{\infty} |\phi_j|\big\}^2  
= \frac{3-\varepsilon}{1-\varepsilon} \|\psi\|_0^2 \|\phi\|_0^2.
\end{equation} Thanks to Remark \ref{rmk:cond}(ii), we have $C_1 \in (0,\infty)$ provided that (SM$_\alpha$) holds for some $\alpha \geq 0$.
For a long memory case, let $K_1, K_2, K_3 \in (0,\infty)$ be the constants associated with $\{\psi_j\}$ and $\{\phi_j\}$ which are defined in Proposition \ref{prop:many-bounds} below and let $K_4 \in (0,\infty)$ and $\{f_k(0)\}_{k=1}^{\infty}$ be defined as in (\ref{eq:K4}) and (\ref{eq:fks}), respectively. 
Let $r \in (1,\infty)$ be chosen such that that $r \sin(\pi d) \in (0,1)$. Using this notation, we define the two constants $C_2, C_3 \in (0,\infty)$ by
\begin{equation} \label{eq:C23} 
C_2 = \frac{2}{1-d} K_1 K_3 \sum_{k=1}^{\infty} f_k(0) \{r\sin(\pi d)\}^k
\quad \text{and} \quad 
C_3 =
K_1K_2K_4 B(d,1-d),
\end{equation} where $B(\cdot, \cdot)$ is the Beta function.

Here are our two versions of a uniform Baxter's inequality.

\begin{theorem}[Uniform Baxter's inequality for short memory time series] \label{thm:unif-BaxterSM}
Suppose $\{X_t\}$ is a short memory time series that satisfies Assumption \ref{assum:SM} for some $\alpha \geq 0$. For $\varepsilon \in (0,1)$, let $C_1 \in (0,\infty)$ be defined as in (\ref{eq:C1}). Then, there exists $N_1 = N_1(\varepsilon) \in \N$ such that 
\begin{equation} \label{eq:unif-Baxter1}
\sum_{k=1}^{n} |\phi_k^m - \phi_{k,n}^m| \leq C_1 \sum_{k=n+1}^{\infty} |\phi_k^m|, \quad m \in \N \quad \text{and} \quad
 n \in \{N_1, N_1+1, ...\}.
\end{equation} The right hand side of above is bounded by
\begin{equation} \label{eq:unif-Baxter2}
\sup_{m\in \N} \sum_{k=n+1}^{\infty} |\phi_k^m| = o(n^{-\alpha})
 \quad \text{as}~~ n \rightarrow \infty.
\end{equation} Furthermore, if $\alpha \geq 1$, then
\begin{equation} \label{eq:unif-Baxter3}
\sum_{m=1}^{\infty} \sum_{k=n+1}^{\infty} |\phi_k^m| = 
o(n^{-\alpha+1}) \quad \text{as}~~ n \rightarrow \infty.
\end{equation}
\end{theorem}

\vspace{0.5em}

Instead of obtaining the bound in terms of the tail sum of $\{\phi_k^m\}_{k=1}^{\infty}$, we will derive the ``exact'' upper bound of $\sum_{k=1}^{n} |\phi_k^m - \phi_{k,n}^m|$ under long memory framework.

\begin{theorem}[Uniform Baxter's inequality for long memory time series] \label{thm:unif-BaxterLM}
Suppose $\{X_t\}$ is a long memory time series that satisfies Assumption \ref{assum:LM} for some $d \in (0,1/2)$. Let $r \in (1,\infty)$ be chosen such that $r \sin(\pi d) \in (0,1)$ and let $C_2,C_3 \in (0,\infty)$ be defined as in (\ref{eq:C23}).
Then, there exists $N_2 = N_2(r) \in \N$ such that 
\begin{equation} \label{eq:unif-BaxterLM1}
\sum_{k=1}^{n} |\phi_k^m - \phi_{k,n}^m| \leq C_2 m^{d} n^{-d}, \quad m \in \N \quad \text{and} \quad
 n \in \{N_2, N_2+1, ...\}.
\end{equation} Moreover, the tail sum of $\{\phi_k^m\}_{k=1}^{\infty}$ is bounded by
\begin{equation} \label{eq:unif-BaxterLM2}
\sum_{k=n}^{\infty} |\phi_k^m| \leq C_3 \left( \frac{m}{n+m}\right)^d, \quad n,m \in \N.
\end{equation} 
\end{theorem}

\subsection{Proof of Theorem \ref{thm:unif-BaxterSM}} \label{appen:pf1}
Given $\varepsilon \in (0,1)$, we choose $N_1 = N_1(\varepsilon) \in \N$ such that $\left \{ \sum_{k=N_1+1}^{\infty} |\phi_k|\right\} \left \{ \sum_{k=1}^{\infty} |\psi_k|\right\} = \varepsilon <1$. This is possible since $ \sum_{k=1}^{\infty} |\psi_k|<\infty$ and $\lim_{n\rightarrow \infty} \sum_{k=n+1}^{\infty} |\phi_k|=0$ under (SM$_\alpha$). By construction, $N_1$ does not depend on $m$. Then, using \cite{SRY}, Lemma B.1 (in Supplementary Material), we have
\begin{equation*} 
\sum_{k=1}^{n} |\phi_k^m - \phi_{k,n}^m| \leq C_{f,0} \sum_{k=n+1}^{\infty} |\phi_k^m|, \quad m \in \N \quad \text{and} \quad
 n \in \{N_1, N_1+1, ...\}.
\end{equation*} Here, $C_{f,0} = \frac{3-\varepsilon}{1-\varepsilon} \|\psi\|_0^2 \|\phi\|_0^2 \in (0,\infty)$, where $\|\psi\|_0$ and $\|\phi\|_0$ are defined as in Remark \ref{rmk:cond}(ii).  Since $C_{f,0}$ does not depend on $m$, (\ref{eq:unif-Baxter1}) holds with $C_{1} = C_{f,0}$.

To show (\ref{eq:unif-Baxter2}), we use (\ref{eq:BLPinfty_ARMA1}) and get
\begin{eqnarray}
\sum_{k=n+1}^{\infty} |\phi_k^m| &\leq& \sum_{k=n+1}^{\infty} \sum_{\ell=0}^{m-1} |\psi_{\ell}| |\phi_{k+m-1-\ell}| 
= \sum_{\ell=0}^{m-1} |\psi_{\ell}| \sum_{k=n+1}^{\infty}  |\phi_{k+m-1-\ell}| \nonumber \\
&\leq& \sum_{\ell=0}^{m-1} |\psi_{\ell}| \sum_{j=n+1}^{\infty}  |\phi_{j}| \leq
\|\psi\|_0 \sum_{j=n+1}^{\infty}  |\phi_{j}|, 
\quad m \in \N \quad \text{and} \quad n \in \N \cup \{0\}.
\label{eq:phikm-SM1}
\end{eqnarray} Therefore, we have
\begin{equation} \label{eq:phikm11}
\sup_{m\in \N} \sum_{k=n+1}^{\infty} |\phi_k^m| \leq \|\psi\|_0 \left\{ \sum_{j=n+1}^{\infty}  |\phi_{j}|\right\}, \quad n \in \N\cup \{0\}.
\end{equation}
Since $\{\phi_j\} \in \mathcal{C}_\alpha^+$, the right hand side of above is bounded by
\begin{equation} \label{eq:bounding}
\sum_{j=n+1}^{\infty}  |\phi_{j}|  \leq (1+n)^{-\alpha} \sum_{j=n+1}^{\infty} (1+j)^\alpha  |\phi_{j}| 
\leq n^{-\alpha} \times o(1) = o(n^{-\alpha}) \quad \text{as}~~ n \rightarrow \infty. 
\end{equation} Substitute above into (\ref{eq:phikm11}) and using that $\|\phi\|_0 <\infty$, we show (\ref{eq:unif-Baxter2}).

Lastly, to show (\ref{eq:unif-Baxter3}), we note that
\begin{eqnarray}
\sum_{m=1}^{\infty} \sum_{k=n+1}^{\infty} |\phi_k^m| &\leq& 
\sum_{m=1}^{\infty} \sum_{\ell=0}^{m-1} |\psi_{\ell}| \sum_{k=n+1}^{\infty}  |\phi_{k+m-1-\ell}| \nonumber \\
&=& \sum_{\ell=0}^{\infty} |\psi_{\ell}|   \sum_{k=n+1}^{\infty}   \sum_{m=\ell+1}^{\infty} |\phi_{k+m-1-\ell}| \quad \text{(Fubini's theorem)} \nonumber \\
&=& \|\psi\|_0 \left\{ \sum_{k=n+1}^{\infty} \sum_{s=0}^{\infty}  |\phi_{k+s}| \right\}
\quad 
\text{(Change of variables: $s=m-1-\ell$)}
\nonumber \\
&\leq& \|\psi\|_0 \left\{ \sum_{u=n+1}^{\infty}  u|\phi_{u}| \right\}.
  \label{eq:phikm-SM2}
\end{eqnarray}
Using a similar technique as in (\ref{eq:bounding}), one can show that $\sum_{u=n+1}^{\infty}  u|\phi_{u}| = o(n^{-\alpha+1})$ as $n\rightarrow \infty$. Therefore, by substituting this bound into (\ref{eq:phikm-SM2}), we show (\ref{eq:unif-Baxter3}).
\hfill $\square$

\subsection{Proof of Theorem \ref{thm:unif-BaxterLM}} \label{appen:pf2}
To prove Theorem \ref{thm:unif-BaxterLM}, we introduce some notation that is used in \cite{p:ino-kas-06}. 
For $n,k \in \N$ and $m, j \in \N \cup \{0\}$, we recursively defined $b_k^m(n,j)$ by
\begin{eqnarray*}
b_1^m(n,j) &=& \phi_j^m = \sum_{k=0}^{m} \psi_k \phi_{j+m-k}, \\
b_{k+1}^m(n,j) &=&\sum_{m_1=0}^{\infty} \phi_{n+1+m_1}^{m} b_{k}^{m_1}(n,j), \quad k \in \N.
\end{eqnarray*} For $m \in \N \cup \{0\}$, $n \in \N$, and $j \in \{1, ...,n\}$, let
\begin{equation} \label{eq:gkm}
g_k^m(n,j) = \begin{cases}
b_k^m(n,j), & k=1, 3, ...\\
b_k^m(n,n+1-j), & k=2, 4, ... \\
\end{cases}.
\end{equation} Then, by \cite{p:ino-kas-06}, Theorem 2.5, the $m$-step ahead finite predictor coefficients can be written as
\begin{equation} \label{eq:phijn-series}
\phi_{j,n}^{m} = \sum_{k=1}^{\infty} g_k^{m-1}(n,j), \quad n,m \in \N \quad \text{and} \quad j \in \{1, 2, ...,n\}.
\end{equation}
By definition, the first term on the right hand side above is
\begin{equation*}
\phi_{j}^{m} = g_1^{m-1}(n,j) = b_j^{m-1}, \quad j,m \in \N.
\end{equation*}

To bound $|g_k^m(n,j)|$, we express $b_k^m(n,j)$ in terms of the absolute convergent series.
Let
\begin{equation*}
\beta_n = \sum_{v=0}^{\infty} \psi_{v} \phi_{n+v}, \quad n \in \N\cup\{0\}.
\end{equation*} Next, we recursively define $\delta_k(n,u,v)$ for $k,n,u,v \in \N \cup\{0\}$ by
\begin{eqnarray*}
\delta_{0}(n,u,v) &=& \delta_{u,v} \\
\text{and} \quad 
\delta_{k+1}(n,u,v) &=& \sum_{w=0}^{\infty} \beta_{n+v+w} \delta_{k}(n,u,w), \quad k\geq 0.
\end{eqnarray*} For example,
$\delta_{3}(n,u,v) = \sum_{s_1=0}^{\infty} \sum_{s_2=0}^{\infty} \beta_{n+u+s_1} \beta_{n+s_1+s_2} \beta_{n+s_2+v}$.
Then, by \cite{p:ino-kas-06}, Theorem 2.8, we have
\begin{equation} \label{eq:thm2.8}
b_k^m(n,j) = \sum_{v=0}^{m} \psi_{m-v} \sum_{u=0}^{\infty} \phi_{j+u} \delta_{k-1}(n+1,u,v), 
\quad n,k \in \N \quad \text{and} \quad m,j \in \N \cup \{0\}.
\end{equation}
Lastly, let $\{f_k(0)\}_{k=1}^{\infty}$ be a sequence of positive numbers that satisfy the following power series expansions:
\begin{equation} \label{eq:fks}
\sum_{k=1}^{\infty} f_{2k}(0) x^{2k} = (\pi^{-1} \arcsin x)^2 \quad \text{and} \quad
\sum_{k=1}^{\infty} f_{2k-1}(0) x^{2k-1} = \pi^{-1} \arcsin x \quad \text{for   } |x|<1.
\end{equation}
 Using the above notation, in the following theorem, we bound $\delta_k(n,u,v)$. 

\begin{theorem} \label{thm:aux1}
Suppose $\{X_t\}$ is a long memory time series that satisfies Assumption \ref{assum:LM} for some $d \in (0,1/2)$. 
For given $r \in (1,\infty)$, there exists $N_2 = N_2(r) \in \N$ such that
\begin{equation} \label{eq:deltak1}
0< \delta_k (n,u,v) \leq \frac{f_k(0) \{r \sin (\pi d)\}^k}{n}, \quad u,v \in \N \cup \{0\}, ~~ k\in\N, ~~ n\geq N_2.
\end{equation} 
\end{theorem}
\noindent \textbf{Proof}. Proof of (\ref{eq:deltak1}) is almost identical to the proof of by \cite{p:ino-kas-06}, Proposition 3.3. Thus, we omit the details. \hfill $\Box$

\vspace{0.5em}

In the remainder of this section, we will provide bounds for the terms involving multistep ahead infinite and finite predictor coefficients. The following proposition is useful for deriving these bounds (the proof is straightforward so we omit the details).
\begin{prop} \label{prop:many-bounds}
Suppose $\{X_t\}$ is a long memory time series that satisfies Assumption \ref{assum:LM} for some $d \in (0,1/2)$. Then, there exist $K_1, K_2, K_3 \in (0,\infty)$ such that
\begin{eqnarray}
\sum_{j=n-1}^{\infty} |\phi_j| &\leq& K_1n^{-d}, \quad n \in \N,  \label{eq:phiu-bound} \\
|\psi_j| &\leq& K_2 (j+1)^{-1+d}, \quad j \in \N \cup \{0\},  \label{eq:K2} \\
\text{and} \quad \sum_{j=1}^{m} |\psi_j| &\leq& K_3 m^d, \quad m \in \N. \label{eq:psiu-bound}
\end{eqnarray}
\end{prop}

Using Theorem \ref{thm:aux1}, we bound $n|\phi_k^m - \phi_{k,n}^m|$ for each $k \in \{1, ...,n\}$.

\begin{theorem}\label{thm:aux2}
Suppose $\{X_t\}$ is a long memory time series that satisfies Assumption \ref{assum:LM} for some $d \in (0,1/2)$. 
Let $N_2 \in \N$ be defined as in (\ref{eq:deltak1}).
 Then, there exists $K >0$ such that for $m \in \N$, $n \in \{N_2, N_2+1, ...\}$, and $k \in \{1, ...,n\}$ such that
\begin{equation} \label{eq:nphi1-bound}
n|\phi_k^m - \phi_{k,n}^m| \leq K \{k \wedge (n+1-k)\}^{-d}m^d,
\end{equation} where $a \wedge b = \min(a,b)$.
\end{theorem} 
\noindent \textbf{Proof}.  
Recall (\ref{eq:phijn-series}) and (\ref{eq:thm2.8}), we have
\begin{eqnarray*}
n(\phi_{k,n}^{m} - \phi_{k}^{m})
&=& n \sum_{j=1}^{\infty} b_{2j+1}^{m-1}(n,k) + n \sum_{j=1}^{\infty} b_{2j}^{m-1}(n,n-k) \\
&=& n \sum_{v=0}^{m-1} \psi_{m-v} \sum_{u=0}^{\infty} \sum_{j=1}^{\infty} \phi_{k+u} \delta_{2j}(n+1,u,v) \\
&&~~+n \sum_{v=0}^{m-1} \psi_{m-v} \sum_{u=0}^{\infty} \sum_{j=1}^{\infty} \phi_{n-k+u} \delta_{2j-1}(n+1,u,v).
\end{eqnarray*} Therefore, for $n,m \in \N$ and $k \in \{1, ...,n\}$, we have
\begin{eqnarray*}
n|\phi_{k,n}^{m} - \phi_{k}^{m}|
&\leq& 
 \sum_{v=0}^{m-1} |\psi_{m-v}| \sum_{u=0}^{\infty}  |\phi_{k+u}| \sum_{j=1}^{\infty} n \delta_{2j}(n+1,u,v) \\
&&~~+ \sum_{v=0}^{m-1} |\psi_{m-v}| \sum_{u=0}^{\infty} |\phi_{n-k+u}|  \sum_{j=1}^{\infty} n\delta_{2j-1}(n+1,u,v)
= I_1(n,m,k) +I_2(n,m,k).
\end{eqnarray*} We bound each term above. Let $r >1$ be chosen so that $0<r\sin(\pi d) <1$. Let $K_1, K_3 \in (0,\infty)$ be defined as in Proposition \ref{prop:many-bounds} and let $\{f_j(0)\}$ be defined as in (\ref{eq:fks}). Let
\begin{equation} \label{eq:K-form}
K= K_1K_3 \sum_{j=1}^{\infty} f_j(0) \{r \sin(\pi d)\}^j.
\end{equation} Then, $K \in (0,\infty)$ due to the exact expressions of the power series (\ref{eq:fks}).
We first bound $I_2(n,m,k)$. For $n\geq N_2$, where $N_2 \in \N$ is as in (\ref{eq:deltak1}),
by using Theorem \ref{thm:aux1} and Proposition \ref{prop:many-bounds},  $I_2(n,m,k)$ is bounded by
\begin{eqnarray}
I_2(n,m,k)&=&\sum_{v=0}^{m-1} |\psi_{m-v}| \sum_{u=0}^{\infty} |\phi_{n-k+u}| \sum_{j=1}^{\infty} n\delta_{2j-1}(n+1,u,v)  \nonumber \\
& \leq & \sum_{v=0}^{m-1} |\psi_{m-v}| \sum_{u=0}^{\infty}  |\phi_{n-k+u}| \sum_{j=1}^{\infty} f_{2j-1}(0) \{r \sin(\pi d)\}^{2j-1} )  \nonumber  \\
&=& \sum_{\ell=1}^{m} |\psi_{\ell}| \sum_{s=n-k}^{\infty} |\phi_s|  \sum_{j=1}^{\infty} f_{2j-1}(0) \{r \sin(\pi d)\}^{2j-1}   )  \nonumber \\
&\leq& K_1 K_3 (n+1-k)^{-d} m^{d} \sum_{j=1}^{\infty} f_{2j-1}(0) \{r \sin(\pi d)\}^{2j-1},
\label{eq:I2}
\end{eqnarray} where $K_1$ and $K_3$ are defined as in (\ref{eq:phiu-bound}) and (\ref{eq:psiu-bound}), respectively. Similarly, for $n \geq N_2$,
$I_1(n,m)$ is bounded by
\begin{eqnarray}
I_1(n,m,k)&=&\sum_{v=0}^{m-1} |\psi_{m-v}| \sum_{u=0}^{\infty}  |\phi_{k+u}| \sum_{j=1}^{\infty} n\delta_{2j}(n+1,u,v) \nonumber \\
& \leq & \sum_{v=0}^{m-1} |\psi_{m-v}| \sum_{s=k}^{\infty}  |\phi_{u}| \sum_{j=1}^{\infty} f_{2j}(0) \{r \sin(\pi d)\}^{2j}
 \nonumber  \\
& \leq & K_1 K_3 k^{-d} m^d \sum_{j=1}^{\infty} f_{2j}(0) \{r \sin(\pi d)\}^{2j}.
\label{eq:I1}
\end{eqnarray} 
Combining (\ref{eq:I2}) and (\ref{eq:I1}), for $m \in \N$, $n \in \{N_2, N_2+1, ...\}$, and $k \in \{1, ...,n\}$, we have
\begin{eqnarray*}
n|\phi_{k,n}^{m} - \phi_{k}^{m}|
&\leq&  I_1(n,m,k) +I_2(n,m,k) \\
&\leq& \bigg[ K_1K_3 \sum_{j=1}^{\infty} f_j(0) \{r \sin(\pi d)\}^j\bigg] \{k \wedge (n+1-k)\}^{-d}m^d \\
&=& K \{k \wedge (n+1-k)\}^{-d}m^d.
\end{eqnarray*} Thus, we get the desired result.
\hfill $\Box$

\vspace{0.5em}

Now, we are ready to prove Theorem \ref{thm:unif-BaxterLM}.

\vspace{0.5em}

\noindent {\bf Proof of Theorem \ref{thm:unif-BaxterLM}}
\hspace{0.3em} We first show (\ref{eq:unif-BaxterLM1}). 
Let $K \in (0,\infty)$ and $N_2 \in \N$ are defined as in (\ref{eq:K-form}) and (\ref{eq:deltak1}), respectively.
By using Theorem \ref{thm:aux2}, we have
\begin{eqnarray*}
\sum_{k=1}^{n} |\phi_{k,n}^{m} - \phi_{k}^{m}| &\leq&
K n^{-1} m^d \sum_{k=1}^{n} \{k \wedge (n+1-k)\}^{-d} 
\leq 2K n^{-1} m^d \sum_{k=1}^{\lfloor(n+1)/2\rfloor} k^{-d} \\
&\leq& 2 (1-d)^{-1}K n^{-1}m^{d} \lfloor(n+1)/2\rfloor^{1-d} \\
&\leq& 2(1-d)^{-1} K m^{d} n^{-d}, \quad m \in \N \quad \text{and} \quad n \in \{N_2, N_2+1, ...\},
\end{eqnarray*} where $\lfloor x\rfloor$ is the greatest greatest integer less than or equal to $x \in \mathbb{R}$. Here, we use
$\sum_{k=1}^{n} k^{-d} \leq (1-d)^{-1} n^{1-d}$, $n \in \N$, in the third inequality and $\lfloor(n+1)/2\rfloor \leq n$ in the last inequality.
 Therefore, (\ref{eq:unif-BaxterLM1}) holds with $C_2 = 2(1-d)^{-1}K \in (0,\infty)$.

To show  (\ref{eq:unif-BaxterLM2}), we use (\ref{eq:BLPinfty_ARMA1}), (\ref{eq:phiu-bound}),
and (\ref{eq:K2}) to get
\begin{eqnarray} 
\sum_{k=n}^{\infty} |\phi_k^m| &\leq&
\sum_{k=n}^{\infty} \sum_{\ell=0}^{m-1} |\psi_\ell| |\phi_{k+m-1-\ell}|
= \sum_{\ell=0}^{m-1}  |\psi_\ell| \sum_{k=n+m-\ell-1}^{\infty}  |\phi_{k}| \nonumber \\
&=&K_1 K_2 \sum_{\ell=0}^{m-1}  \frac{1}{(1+\ell)^{1-d}} \frac{1}{(n+m-\ell)^d} \nonumber \\
&=& K_1 K_2 \sum_{\ell=1}^{m}  \frac{1}{\ell^{1-d}} \frac{1}{(n+m-\ell+1)^d}, \quad n,m \in \N.
\label{eq:pkm22}
\end{eqnarray} Since $f(\ell) = (m-\ell+1)/(n+m-\ell+1)$ is a decreasing function of $\ell$, we have
\begin{equation*}
\frac{(m-\ell+1)^d}{(n+m-\ell+1)^d} = \{f(\ell)\}^d \leq \{f(1)\}^d = \frac{m^d}{(n+m)^d}, \quad n,m \in \N, ~~ \ell \in \{1, ...,m\}.
\end{equation*} Substitute the above into (\ref{eq:pkm22}), 
\begin{eqnarray} 
\sum_{k=n}^{\infty} |\phi_k^m| &\leq&
K_1 K_2  \left( \frac{m}{n+m}\right)^d \sum_{\ell=1}^{m}  \frac{1}{\ell^{1-d}} \frac{1}{(m-\ell+1)^d} \nonumber \\
&=& K_1 K_2 \left( \frac{m}{n+m}\right)^d  \frac{1}{m}\sum_{\ell=1}^{m}  \frac{1}{(\ell/m)^{1-d}} \frac{1}{\{(1+1/m)-\ell/m\}^d}
, \quad n,m \in \N.
\label{eq:pkm33}
\end{eqnarray} Finally, since $\lim_{m\rightarrow \infty} \frac{1}{m}\sum_{\ell=1}^{m}  \frac{1}{(\ell/m)^{1-d}} \frac{1}{\{(1+1/m)-\ell/m\}^d} = \int_{0}^{1} x^{d-1}(1-x)^{-d} dx = B(d,1-d)$, where $B(\cdot, \cdot)$ is the Beta function,
there exists $K_4 \in (0,\infty)$ such that for $m\in \N$,
\begin{equation} \label{eq:K4}
\frac{1}{m}\sum_{\ell=1}^{m}  \frac{1}{(\ell/m)^{1-d}} \frac{1}{\{(1+1/m)-\ell/m\}^d} \leq K_4 \int_{0}^{1} x^{d-1} (1-x)^{-d} dy = K_4 B(d,1-d).
\end{equation}
Substitute the above into (\ref{eq:pkm33}), we have
\begin{equation*}
\sum_{k=n}^{\infty} |\phi_k^m| \leq K_1 K_2 K_4 B(d,1-d) \left( \frac{m}{n+m}\right)^d, \quad n,m\in \N.
\end{equation*} Therefore,  (\ref{eq:unif-BaxterLM2}) holds with $C_3 =K_1 K_2 K_4 B(d,1-d) \in (0,\infty)$.
All together, we prove the theorem.
\hfill $\Box$

\section{Proof of the main results} \label{sec:proof-main}

\subsection{Proof of Theorem \ref{thm:expression}}
For $n \in \N$, let $Y_n=  H(B) X_n = \sum_{k=-\infty}^{\infty} h_k X_{n-k}$ be a target random variable. 
Since $P_{(-\infty, n]} X_{n-m} = X_{n-m}$ for $n \in \N$ and $m \in \N \cup \{0\}$, we have
\begin{eqnarray*}
P_{(-\infty, n]} Y_n
&=& P_{(-\infty, n]} \sum_{k=-\infty}^{\infty}   h_k X_{n-k} 
= \sum_{k=-\infty}^{-1}  h_k P_{(-\infty, n]}  X_{n-k} + \sum_{k=0}^{\infty}h_k X_{n-k} \\
&=&\sum_{m=1}^{\infty}  h_{-m} P_{(-\infty, n]} X_{n+m} + \sum_{k=0}^{\infty} h_k X_{n-k} \quad \text{(Change of variables: $m=-k$)} \\
&=&\sum_{m=1}^{\infty} h_{-m} \sum_{k=1}^{\infty} \phi_{k}^m  X_{n+1-k} + \sum_{k=0}^{\infty} h_k X_{n-k} 
\quad \text{(From (\ref{eq:m-step1}))} \\
&=& \sum_{k=1}^{\infty} \left( \sum_{m=1}^{\infty} h_{-m} \phi_{k}^m \right) X_{n+1-k} + \sum_{k=0}^{\infty}h_{k} X_{n-k},
\quad n \in\N.
\end{eqnarray*} The absolute convergence of $\sum_{m=1}^{\infty} h_{-m} \phi_{k}^m$ is due to Lemma \ref{lemma:abs-conv}.
Comparing the coefficients above with $\{\hat{h}_k\}_{k=1}^{\infty}$ in (\ref{eq:Yhat}), we show (\ref{eq:FIRcoeff1}).

To show (\ref{eq:FIRcoeff2}), we first note that
\begin{equation*}
P_{[1, n]} X_{n-m} = X_{n-m}, \quad n \in \N \quad \text{and} \quad m \in \{0,1, ...,n-1\}.
\end{equation*} Therefore,
\begin{eqnarray*}
P_{[1, n]} Y_n
&=& P_{[1, n]} \sum_{k=-\infty}^{\infty}   h_k X_{n-k} \\
&=& \sum_{k=-\infty}^{-1}  h_k P_{[1, n]}  X_{n-k} + \sum_{k=0}^{n-1} h_k X_{n-k} + \sum_{k=n}^{\infty}  h_k P_{[1, n]}  X_{n-k} \\
&=& \sum_{m=1}^{\infty} h_{-m} P_{[1, n]}  X_{n+m} + \sum_{k=0}^{n-1} h_k X_{n-k} 
+ \sum_{m=1}^{\infty}  h_{n+m-1}  P_{[1, n]}  X_{1-m}.
\end{eqnarray*} By using (\ref{eq:m-step2}), the first term above is
\begin{equation*}
\sum_{m=1}^{\infty} h_{-m} P_{[1, n]}  X_{n+m} =
\sum_{m=1}^{\infty} h_{-m} \sum_{k=1}^{n} \phi_{k,n}^m  X_{n+1-k} 
=
\sum_{k=1}^{n} \left( \sum_{m=1}^{\infty} h_{-m} \phi_{k,n}^m \right) X_{n+1-k}.
\end{equation*} Using the reflexive property, that is,
$P_{[1, n]}  X_{1-m} = \sum_{k=1}^{n} \phi_{n+1-k,n}^{m+1}   X_{n+1-k}$ for $n,m \in \N$, 
the third term above is
\begin{eqnarray*}
\sum_{m=1}^{\infty}  h_{n+m-1}  P_{[1, n]}  X_{1-m} &=&
\sum_{m=1}^{\infty} h_{n+m-1}  \sum_{k=1}^{n} \phi_{n+1-k,n}^{m}  X_{n+1-k} \\
&=&
\sum_{k=1}^{n} \left( \sum_{m=1}^{\infty} h_{n+m-1} \phi_{n+1-k,n}^{m} \right) X_{n+1-k}.
\end{eqnarray*} All together, we have
\begin{equation*}
P_{(-\infty, n]} Y_n 
= \sum_{k=1}^{n} \left( \sum_{m=1}^{\infty} h_{-m} \phi_{k,n}^m \right) X_{n+1-k}  +  \sum_{k=1}^{n} h_{k-1} X_{n+1-k} 
+\sum_{k=1}^{n} \left( \sum_{m=1}^{\infty} h_{n+m-1} \phi_{n+1-k,n}^{m} \right) X_{n+1-k}.
\end{equation*} The absolute convergence of $\sum_{m=1}^{\infty} h_{-m} \phi_{k,n}^m$ and $\sum_{m=1}^{\infty} h_{n+m-1} \phi_{n+1-k,n}^{m}$ are due to Lemma \ref{lemma:abs-conv2}.
Comparing the coefficients above with $\{\hat{h}_{k,n}\}_{k=1}^{n}$ in (\ref{eq:Yhatn}), we show (\ref{eq:FIRcoeff2}).

Lastly, (\ref{eq:hatH1}) and (\ref{eq:hatH2}) are due to (\ref{eq:abs-hatH1}) and (\ref{eq:abs-hatH2}) in the Appendix, respectively. All together, we prove the theorem.
\hfill $\square$

\subsection{Proof of Theorem \ref{thm:L1-short}}

Let $C_1 \in (0,\infty)$ and $N_1 \in \N$ are defined as in (\ref{eq:C1}) and (\ref{eq:unif-Baxter1}), respectively. Then, from (\ref{eq:hatHdiff}) and (\ref{eq:unif-Baxter1}), we have
\begin{eqnarray*}
\sum_{k=1}^{n} |\hat{h}_k - \hat{h}_{k,n}| 
&\leq&  \sum_{m=1}^{\infty}  | h_{-m}|  \sum_{k=1}^{n} | \phi_{k}^m - \phi_{k,n}^m| +
 \sum_{m=1}^{\infty} \sum_{k=1}^{n}  |h_{n+m-1}|   |\phi_{n+1-k,n}^{m}|
\\
&\leq& C_1 \sum_{m=1}^{\infty}  | h_{-m}|  \sum_{k=n+1}^{\infty} | \phi_{k}^m|
+ \sum_{m=1}^{\infty} \sum_{k=1}^{n}  |h_{n+m-1}|   |\phi_{n+1-k,n}^{m}|
, \quad n \in \{N_1, N_1+1, ...\}.
\end{eqnarray*} If $\{h_k\} \in \ell_1$, then by using (\ref{eq:phikm11}) and (\ref{eq:tail-SM1}) (of Lemma \ref{lemma:abs-conv2}), we have
\begin{equation*}
\sum_{m=1}^{\infty}  | h_{-m}|  \sum_{k=n+1}^{\infty} | \phi_{k}^m|
\leq \left\{ \sup_{m\in\N} \sum_{k=n+1}^{\infty} | \phi_{k}^m|\right\} \sum_{m=1}^{\infty}  | h_{-m}|
 = 
\|\psi\|_0 \sum_{m=0}^{\infty} |h_{-m}| \times
o(n^{-\alpha})
\end{equation*} and
\begin{equation*}
\sum_{m=1}^{\infty} \sum_{k=1}^{n}  |h_{n+m-1}|   |\phi_{n+1-k,n}^{m}|
 \leq \left\{  \|\psi\|_0 \|\phi\|_0+ o(n^{-\alpha}) \right\} \times \sum_{k=n}^{\infty} |h_k|
 \quad \text{as}~~ n \rightarrow \infty.
\end{equation*} Combining, we have
\begin{eqnarray}
\sum_{k=1}^{n} |\hat{h}_k - \hat{h}_{k,n}|  &\leq&
C_1 \|\psi\|_0 \sum_{m=0}^{\infty} |h_{-m}| \times o(n^{-\alpha}) + 
\left\{  \|\psi\|_0 \|\phi\|_0+ o(n^{-\alpha}) \right\} \times  \sum_{k=n}^{\infty} |h_k|  \nonumber \\
&=& o(n^{-\alpha}) + O\bigg( \sum_{k=n}^{\infty} |h_k| \bigg) \quad \text{as}~~ n \rightarrow \infty.
\label{eq:SM1-exactbound}
\end{eqnarray}
Thus, we show (\ref{eq:short-baxter11}).

To show (\ref{eq:short-baxter22}), by using (\ref{eq:phikm-SM2}) and (\ref{eq:tail-SM2}) (in the Appendix), we have
\begin{equation*}
\sum_{m=1}^{\infty}  | h_{-m}|  \sum_{k=n+1}^{\infty} | \phi_{k}^m|
\leq \left\{ \sup_{m\in\N} |h_{-m}| \right\} \sum_{m=1}^{\infty}  \sum_{k=n+1}^{\infty} | \phi_{k}^m|
 = \sup_{m\in\N} |h_{-m}| \cdot \|\psi\|_0 \times o(n^{-\alpha+1})
\end{equation*} 
and
\begin{equation*}
\sum_{m=1}^{\infty} \sum_{k=1}^{n}  |h_{n+m-1}|   |\phi_{n+1-k,n}^{m}|
 = \big\{ \|\psi\|_0 \|\phi\|_1  + o(n^{-\alpha+1}) \big\} \times \sup_{k \geq n} |h_k| \quad \text{as}~~ n \rightarrow \infty.
\end{equation*} Combining, 
\begin{eqnarray}
\sum_{k=1}^{n} |\hat{h}_k - \hat{h}_{k,n}|  &\leq&
C_1 \|\psi\|_0 \sup_{m\in\N} |h_{-m}| \times o(n^{-\alpha+1}) +   
\big\{ \|\psi\|_0 \|\phi\|_1  + o(n^{-\alpha+1}) \big\} \times \sup_{k \geq n} |h_k|  \nonumber \\
&=& o(n^{-\alpha+1}) + O\big( \sup_{k \geq n} |h_k| \big) \quad \text{as}~~ n \rightarrow \infty.
\label{eq:SM2-exactbound}
\end{eqnarray}
Thus, we show (\ref{eq:short-baxter22}). All together, we get the desired results.
\hfill $\Box$

\subsection{Proof of Theorem \ref{thm:L1-long}}

 Let $C_2,C_3 \in (0,\infty)$ be defined as in (\ref{eq:C23}) and let $N_2 \in \N$ be defined as in (\ref{eq:unif-BaxterLM1}). Then, from (\ref{eq:hatHdiff}) and (\ref{eq:unif-BaxterLM1}), we have
\begin{eqnarray*}
\sum_{k=1}^{n} |\hat{h}_k - \hat{h}_{k,n}| 
&\leq&  \sum_{m=1}^{\infty}  | h_{-m}|  \sum_{k=1}^{n} | \phi_{k}^m - \phi_{k,n}^m| +
 \sum_{m=1}^{\infty} \sum_{k=1}^{n}  |h_{n+m-1}|   |\phi_{n+1-k,n}^{m}| \\
&\leq& C_2 n^{-d} \sum_{m=1}^{\infty}  | h_{-m}| m^{d}
+  \sum_{m=1}^{\infty} \sum_{k=1}^{n}  |h_{n+m-1}|   |\phi_{n+1-k,n}^{m}|, \quad n \in \{N_2, N_2+1, ...\}.
\end{eqnarray*} Using that $(h_{-1}, h_{-2}, ....) \in \mathcal{C}_d^{-}$ and (\ref{eq:tail-LM}) (in the Appendix), the above is bounded by
\begin{equation}
\sum_{k=1}^{n} |\hat{h}_k - \hat{h}_{k,n}|  \leq
C_2 \sum_{m=1}^{\infty}  | h_{-m}| m^{d} \cdot n^{-d} +   
\{C_2+C_3\} \times o(n^{-d}) = O(n^{-d}), \quad  n \rightarrow \infty.
\label{eq:LM-exactbound}
\end{equation}
Thus, we get the desired result.
\hfill $\Box$

\subsection{Proof of Theorem \ref{thm:sigma-SM}}
A proof is straightforward by using (\ref{eq:sigmatilde-bound}), (\ref{eq:sigma-bound}), Theorem \ref{thm:L1-short}, and Lemma \ref{lemma:abs-conv3} in the Appendix (under (H-SM1) and (H-SM2)).
We note that by using precise expressions in (\ref{eq:SM1-exactbound}) and (\ref{eq:tail-SM1})  (resp., (\ref{eq:SM2-exactbound}) and (\ref{eq:tail-SM2})), one can easily derive the constant factors that are multiplied by the asymptotic notation in (\ref{eq:sigma-SM1}) (resp., (\ref{eq:sigma-SM2})).
\hfill $\Box$

\subsection{Proof of Theorem \ref{thm:sigma-LM}}

A proof is straightforward by using (\ref{eq:sigmatilde-bound}), (\ref{eq:sigma-bound}), Theorem \ref{thm:L1-long}, and Lemma \ref{lemma:abs-conv3} in the Appendix (under (H-LM)). 
We note that by using precise expressions in (\ref{eq:LM-exactbound}) and (\ref{eq:tail-LM}), one can easily derive the constant factors that are multiplied by the asymptotic notation in (\ref{eq:sigma-LM}).
\hfill $\Box$

\section{A multivariate extension} \label{sec:multi}

One possible extension of our results is in the multivariate framework. Before introducing the multivariate analogue of the approximation problem, we first fix terms. Let $\{\underline{X}_t\} = \{\underline{X}_t\}_{t \in \Z}$ be a $q$-variate centered second-order stationary process with components $\underline{X}_t = (X_t^1, ..., X_t^q)^\top$. The corresponding $q\times q$ matrix-valued autorcovariate and spectral density functions are denoted by $\Gamma(k) = \cov(X_k, X_0)$ and $F(e^{i\omega}) = \sum_{k\in \Z} \Gamma(k)e^{ik\omega}$, respectively. For the notational convenience, we will use the same notation as in Section \ref{sec:prelim} (under univariate setting) to define the Hilbert space associated with $\{\underline{X}_t\}$.
Let $V$ be the real Hilbert space spanned by $\{X_t^k: t \in \Z, ~ 1\leq k \leq q\}$,
 and let $V^q$ be the space of $q$-variate random variables where each element belongs to $V$. Then, $V^q$ is also a Hilbert space with the norm $\|(x_1, ...,x_q)^\top \|_{V^q} = \{\sum_{k=1}^{q} \|x_k\|_{V}^2\}^{1/2}$. For $K \subset \Z$ and $\underline{x} = (x_1, ... ,x_q)^\top \in V^q$, let $P_K \underline{x} = (P_K x_1, ..., P_K x_q)^\top$, where $P_K$ is a projection of $V$ onto the closed subspace $V_K =\cspn \{X_t^k:  t \in K,~ 1\leq k \leq q \}$.

Given an $\R^{q\times q}$-valued filter coefficients $\{H_k\}_{k\in \Z}$, let $\underline{Y}_n = \sum_{k \in \Z} H_k \underline{X}_{n-k} \in \R^{q}$ ($n \in \N$) be a multivariate target random variable. For important examples of the (multivariate) filter coefficients and the corresponding target random variables, we refer readers to \cite{p:mce-wil-20}, Section 3. Analogous to the univariate definitions of $\widehat{Y}_n$, $\widehat{Y}_n^{(n)}$, and $\widetilde{\underline{Y}}_n^{(n)}$ in (\ref{eq:Yhat}), (\ref{eq:Yhatn}), and (\ref{eq:Ytilde}), respectively, for each $n \in \N$, we define
\begin{eqnarray*} 
&& \widehat{\underline{Y}}_n = P_{(-\infty,n]} \underline{Y}_n = \sum_{k=1}^{\infty} \widehat{H}_k \underline{X}_{n+1-k}, \quad \widehat{Y}_n^{(n)} = P_{[1,n]} \underline{Y}_n = \sum_{k=1}^{n} \widehat{H}_{k,n} \underline{X}_{n+1-k}, \\
&& \text{and} \quad
\widetilde{\underline{Y}}_n^{(n)} = \sum_{k=1}^{n} \widehat{H}_k \underline{X}_{n+1-k},
\end{eqnarray*} 
 where $\{\widehat{H}_k\}_{k=1}^{\infty}$ and $\{\widehat{H}_{k,n}\}_{k=1}^{n}$ are the best (linear) predictor coefficients of $\underline{Y}_n$ given $V_{(-\infty,n]}$ and $V_{[1,n]}$, respectively. Then, we are interested in the rate of convergence to zero of the following three quantities: $\sum_{k=1}^{n}\|\widehat{H}_k - \widehat{H}_{k,n}\|_{\text{spec}}$, 
$\widetilde{\delta}_n = \| \widehat{\underline{Y}}_n^{(n)} - \widetilde{\underline{Y}}_n^{(n)}\|_{V^q}$, and
$\delta_n= \| \widehat{\underline{Y}}_n^{(n)} - \widehat{\underline{Y}}_n\|_{V^q}$. Here $\|\cdot\|_{\text{spec}}$ denotes the spectral norm.
We believe that all the asymptotic results in Section \ref{sec:main-results} and the inequalities in Section \ref{sec:Baxter-main} should transfer over to the multivariate analogues. However, there are two subtle differences compared to the univariate case, which we will outline below.

\vspace{0.5em}

\noindent $\bullet$ \textit{Wold-type representations and the assumptions}. The Wold-type infinite-order MA and AR representations of $\{X_t\}$ in Section \ref{sec:prelim} are the building blocks of our theoretical development in the univariate framework.
For the multivariate case, we can obtain similar representations. Specifically, under the conditions stated in \cite{b:roz-67}, Chapter II (see also \cite{IKP2}, pages 1208--1209), $\{\underline{X}_t\}$ admits the following \textit{forward} MA and AR representations:
\begin{equation} \label{eq:multi-Wold1}
\underline{X}_t = \underline{\varepsilon}_{t} + \sum_{j=1}^{\infty} \Psi_j \underline{\varepsilon}_{t-j} \quad \text{and} \quad
\underline{X}_{t} - \sum_{j=1}^{\infty} \Phi_j \underline{X}_{t-j} = \underline{\varepsilon}_{t}, \quad t\in\Z,
\end{equation} where $\{\underline{\varepsilon}_{t}\}$ is a $q$-variate white noise process and
$\{\Psi_j\}_{j=0}^{\infty}$ (setting $\Psi_0 = I_q$, the identity matrix of order $q \in \N$) and $\{\Phi_j\}_{j=0}^{\infty}$ (setting $\Phi_0 = -I_q$)  are the $\R^{q \times q}$-valued infinite-order forward MA and AR coefficients of $\{\underline{X}_t\}$, respectively. 

 Next, we introduce the analogous class of multivariate short and long memory processes. For the short memory case, we denote the condition (SM$_\alpha^\prime$) for some $\alpha \in [0,\infty)$ by
\begin{equation} \tag{SM$_\alpha^\prime$}
\{\Psi_j\}, \{\widetilde{\Psi}_j\}, \{\Phi_j\}, \{\widetilde{\Phi}_j\}  \in \mathcal{C}_\alpha^{q\times q}(\N \cup \{0\}),
\end{equation} where 
$\mathcal{C}_\alpha^{q\times q}(\N \cup \{0\}) = \{(A_k \in \R^{q\times q}: k\geq 0): \sum_{k=0}^{\infty} (1+|k|)^\alpha \|A_k\|_{\text{spec}} < \infty \}$. Here, $ \{\widetilde{\Psi}_j\}$ and $\{\widetilde{\Phi}_j\}$ are the infinite-order \textit{backward} MA and AR coefficients of $\{\underline{X}_t\}$, respectively (see \cite{IKP2}, pages 1208--1209 for the precise definition). For the long-memory case, we will focus on the spectral density matrices of form
\begin{equation}  \tag{LM$_d^\prime$}
F(e^{i\omega}) = |1-e^{i\omega}|^{-2d} g(e^{i\omega}) g(e^{i\omega})^* \quad \text{for some} \quad d\in (0,1/2),
\end{equation} where $g(z)$ is a rational function (on $z \in \C$) that satisfies the condition (C) in \cite{IKP2} (page 1221). Heuristically, $F(\cdot)$ corresponds to the spectral density matrix of a $q$-variate ARFIMA process with a common differencing parameter $d \in (0,1/2)$. Then, similarly to the condition (LM$_d$) for univariate case, $\|\Psi_j\|_{\text{spec}}$ and $\|\widetilde{\Psi}_j\|_{\text{spec}}$ decreasing at a rate $O(j^{-1+d})$ and $\|\Phi_j\|_{\text{spec}}$ and $\|\widetilde{\Phi}_j\|_{\text{spec}}$ decreasing at a rate $O(j^{-1-d})$ as $j \rightarrow \infty$ (cf. \cite{I00}, Lemma 2.2). 

Lastly, we define one of the three analogous conditions for the filter coefficients $\{H_k\}$ in the multivariate setting. An analogous version of (H-SM1) in Assumption \ref{assum:H} is 
\begin{equation} \tag{H-SM1$^\prime$}
\text{$\{X_t\}$ satisfies (SM$_\alpha^\prime$) for some $\alpha \geq 0$ and } \{H_k\} \in \ell_1^{q\times q}(\Z),
\end{equation} where $\ell_p^{q \times q}(\Z) = \{(A_k \in \R^{q\times q}: k \in \Z): \sum_{k \in \Z} \|A_k\|_{\text{spec}}^{p} < \infty \}$, $p \in [1,\infty)$. Analogous versions of (H-SM2) and (H-LM) in Assumption \ref{assum:H} for the multivariate setting can be defined in a same manner.

\vspace{0.5em}

\noindent $\bullet$ \textit{Uniform Baxter's inequality}. A key ingredient of the proof of the asymptotic rates in Section \ref{sec:main-results} is a uniform Baxter's inequality under both short memory (Theorem \ref{thm:unif-BaxterSM}) and long memory (Theorem \ref{thm:unif-BaxterLM}) time series. We thus need to prove the analogous uniform Baxter's inequalities for multivariate time series. Here, we will only focus on the \textit{forward} predictor coefficients (coefficients that appear in $P_{(-\infty,n]}\underline{X}_{n+m}$ and  $P_{[1,n]}\underline{X}_{n+m}$) and the \textit{backward} predictor coefficients (coefficients that appear in $P_{[0,\infty)}\underline{X}_{-m}$ and  $P_{[0,n-1]}\underline{X}_{-m}$) can be treated in the same way. Under (SM$_\alpha^\prime$), \cite{p:che-93} generalized the Baxter's inequality for the one-step forward predictor coefficients (cf. equation (\ref{eq:Baxter})) to the multivariate case. Also, in \cite{IKP2}, they provided a series expansion of the forward (and backward) one-step finite predictor coefficients (cf. equation (\ref{eq:phijn-series})) and rigorously proved the analogous version of Baxter's inequality under (LM$_d^\prime$) for some $d \in (0,1/2)$. The theoretical properties derived in the aforementioned articles together with the techniques developed in Sections \ref{appen:pf1} and \ref{appen:pf2} (also in \cite{SRY}, Lemma A.2) may be useful in showing the analogous uniform Baxter's inequality for multivariate short and long memory time series. This will be investigated in future research.

\section*{Acknowledgments}
The author's research was supported by Taiwan's National Science and Technology Council (grant 110-2118-M-001-014-MY3). 
The author is grateful to Professor Akihiko Inoue for fruitful discussions.
The author also wishes to thank the two anonymous referees and editors for their valuable comments and corrections,
 which have greatly improved the article in all aspects.

\begin{appendix}

\section{Technical lemmas} \label{appen:A}

The following three auxiliary lemmas are required in the proof of Section \ref{sec:main-results}. Lemma \ref{lemma:abs-conv} is used to show (\ref{eq:hatH1}) and (\ref{eq:hatH2}). Lemma \ref{lemma:abs-conv2} is used to bound the second term of an expansion of $\sum_{k=1}^{n} |\hat{h}_k - \hat{h}_{k,n}|$ in (\ref{eq:hatHdiff}).  Lemma \ref{lemma:abs-conv3} is used to bound the second term of an expansion of $\sigma_n$ in (\ref{eq:sigma-bound}).

\begin{lemma} \label{lemma:abs-conv}
Suppose $\{X_t\}$ satisfies either Assumption \ref{assum:SM} or \ref{assum:LM} and the filter coefficients $\{h_k\}$ satisfy Assumption \ref{assum:H}. Then,
under (H-SM1) or (H-LM),
\begin{equation} \label{eq:abs-hatH1}
\sum_{k=1}^{\infty} \left\{ 
\sum_{m=1}^{\infty} |h_{-m}|  |\phi_{k}^m|  + |h_{k-1}|
\right\} <\infty
\end{equation}  and under (H-SM2),
\begin{equation}  \label{eq:abs-hatH2}
\sum_{k=1}^{\infty} \left\{ 
\sum_{m=1}^{\infty} |h_{-m}|  |\phi_{k}^m|  + |h_{k-1}|
\right\}^2 <\infty.
\end{equation}
\end{lemma}
\noindent \textbf{Proof}. We first show (\ref{eq:abs-hatH1}).  Since $(h_0,h_1, ...) \in \ell_1^+$ under (H-SM1) and (H-LM), it is sufficient to show
\begin{equation} \label{eq:hmfubini}
\sum_{m=1}^{\infty}  |h_{-m}| \sum_{k=1}^{\infty} |\phi_k^m| < \infty.
\end{equation} 
First, we assume (H-SM1). Then, by using (\ref{eq:phikm-SM1}), we have
\begin{equation*}
\sup_{m \in \N} \sum_{k=1}^{\infty} |\phi_k^m|  \leq
 \sum_{\ell=0}^{\infty} |\psi_{\ell}| \sum_{j=0}^{\infty}  |\phi_{j}| <\infty.
\end{equation*} 
Under (H-LM), by using (\ref{eq:unif-BaxterLM2}), 
\begin{equation*}
\sup_{m \in \N} \sum_{k=1}^{\infty} |\phi_k^m|  \leq C_3 \sup_{m \in \N} \left(\frac{m}{1+m} \right)^d = C_3<\infty.
\end{equation*}
therefore, since $(h_{-1}, h_{-2}, ...) \in \ell_1^-$, under either (H-SM1) or (H-LM), we have
\begin{equation*}
\sum_{m=1}^{\infty}  |h_{-m}| \sum_{k=1}^{\infty} |\phi_k^m| \leq
\sup_{m \in \N} \sum_{k=1}^{\infty} |\phi_k^m| \cdot \sum_{m=1}^{\infty}  |h_{-m}| <\infty.
\end{equation*} Lastly, by Fubini's theorem, (\ref{eq:abs-hatH1}) is followed.

Next, we will show (\ref{eq:abs-hatH2}). By iteratively applying the Cauchy-Schwarz inequality, we have
\begin{eqnarray} 
\sum_{k=1}^{\infty} \left\{ 
\sum_{m=1}^{\infty} |h_{-m}|  |\phi_{k}^m|  + |h_{k-1}|
\right\}^2 &\leq& 2 \sum_{k=1}^{\infty} \left\{ \sum_{m=1}^{\infty} |h_{-m}|  |\phi_{k}^m| \right\}^2
+ 2 \sum_{k=1}^{\infty} |h_{k-1}|^2 \nonumber  \\
&\leq&
 2 \sum_{m=1}^{\infty} |h_{-m}|^2 \cdot   \sum_{k=1}^{\infty} \left\{ \sum_{m=1}^{\infty} |\phi_{k}^m|^2 \right\}
+ 2 \sum_{k=1}^{\infty} |h_{k-1}|^2 \nonumber  \\
&\leq&
 2  \sum_{m=1}^{\infty} |h_{-m}|^2 \cdot  \left\{ \sum_{k=1}^{\infty} \sum_{m=1}^{\infty} |\phi_{k}^m| \right\}^2
+ 2 \sum_{k=1}^{\infty} |h_{k-1}|^2.
\label{eq:hsq11}
\end{eqnarray}
Under (H-SM2), by (\ref{eq:phikm-SM2}), we have
\begin{equation*}
\sum_{m=1}^{\infty} \sum_{k=1}^{\infty} |\phi_k^m| \leq 
\|\psi\|_0 \cdot \sum_{u=1}^{\infty}  u|\phi_{u}| <\infty.
\end{equation*} Therefore, combining the above two inequalities, we have
\begin{equation*} 
\sum_{k=1}^{\infty} \left\{ 
\sum_{m=1}^{\infty} |h_{-m}|  |\phi_{k}^m|  + |h_{k-1}|
\right\}^2
\leq 2 \|\psi\|_0^2 \left\{ \sum_{m=1}^{\infty} |h_{-m}|^2\right\} \left\{ \sum_{u=1}^{\infty}  u|\phi_{u}| \right\}^2
+ 2 \sum_{k=1}^{\infty} |h_{k-1}|^2 <\infty.
\end{equation*}
Thus, we get the desired result.
\hfill $\Box$

\begin{lemma} \label{lemma:abs-conv2}
Suppose $\{X_t\}$ satisfies either Assumption \ref{assum:SM} or \ref{assum:LM} and the filter coefficients $\{h_k\}$ satisfy Assumption \ref{assum:H}. Then, there exists $N\in \N$ such that 
\begin{equation} \label{eq:hmbound11}
\sum_{m=1}^{\infty} |h_{-m}|  |\phi_{k,n}^m| <\infty, \quad n \in \{N,N+1, ...\} \quad \text{and} \quad k \in \{1,2, ...,n\}.
\end{equation} Moreover,
\begin{equation} \label{eq:hmbound22}
\sum_{m=1}^{\infty} \sum_{k=1}^{n}  |h_{n+m-1}|   |\phi_{n+1-k,n}^{m}| 
=
\begin{cases}
O(\sum_{k=n}^{\infty} |h_{k}|), & \text{under (H-SM1)} \\
O(\sup_{k \geq n} |h_k|) , 
& \text{under (H-SM2)} \\
o(n^{-d}), & \text{under (H-LM)} 
\end{cases} \quad \text{as}~~ n\rightarrow \infty.
\end{equation} 
\end{lemma}
\noindent \textbf{Proof}.
To show (\ref{eq:hmbound11}), we use a triangular inequality 
\begin{eqnarray*}
\sum_{m=1}^{\infty} |h_{-m}| |\phi_{k,n}^m|
&\leq& \sum_{m=1}^{\infty} |h_{-m}| |\phi_{k}^m|  + \sum_{m=1}^{\infty} |h_{-m}| |\phi_{k}^m - \phi_{k,n}^m| \\
&\leq& \sum_{m=1}^{\infty} |h_{-m}| |\phi_{k}^m|  + \sum_{m=1}^{\infty} |h_{-m}| \sum_{k=1}^{n} |\phi_{k}^m - \phi_{k,n}^m|.
\end{eqnarray*} The first term above is finite due to Lemma \ref{lemma:abs-conv}. To show the absolute convergence of the second term,
let $N = \max \{N_1,N_2\} \in \N$, where $N_1$ an $N_2$ are defined as in (\ref{eq:unif-Baxter1}) and (\ref{eq:unif-BaxterLM1}), respectively. Then,
 under (H-SM1), by using (\ref{eq:unif-Baxter1}) and
(\ref{eq:unif-Baxter2}), for $n \geq N$, we have
\begin{equation} \label{eq:HSM1-bound}
\sum_{m=1}^{\infty} |h_{-m}| \sum_{k=1}^{n} |\phi_{k}^m - \phi_{k,n}^m| \leq
\sum_{m=1}^{\infty} |h_{-m}| \sum_{k=n+1}^{\infty} |\phi_{k}^m| 
\leq \sup_{m\in\N} \sum_{k=n+1}^{\infty} |\phi_{k}^m| \sum_{m=1}^{\infty} |h_{-m}| <\infty,
\end{equation} Under (H-SM2), by using (\ref{eq:unif-Baxter1}) and
(\ref{eq:unif-Baxter3}), we have
\begin{equation} \label{eq:HSM2-bound}
\sum_{m=1}^{\infty} |h_{-m}| \sum_{k=1}^{n} |\phi_{k}^m - \phi_{k,n}^m| \leq
\sup_{m\in\N}  |h_{-m}|  \sum_{m=1}^{\infty}  \sum_{k=n+1}^{\infty} |\phi_{k}^m| <\infty,~~ n \geq N.
\end{equation} Finally, under (H-LM), by using (\ref{eq:unif-BaxterLM1}),
\begin{equation} \label{eq:HLM-bound}
\sum_{m=1}^{\infty} |h_{-m}| \sum_{k=1}^{n} |\phi_{k}^m - \phi_{k,n}^m| \leq
C_2 n^{-d} \sum_{m=1}^{\infty} |h_{-m}| m^{d} <\infty,~~ n \geq N.
\end{equation} Therefore, combining (\ref{eq:HSM1-bound}), (\ref{eq:HSM2-bound}), and (\ref{eq:HLM-bound}), we show (\ref{eq:hmbound11}).

To show (\ref{eq:hmbound22}), we note that
\begin{eqnarray*}
\sum_{m=1}^{\infty} \sum_{k=1}^{n}  |h_{n+m-1}|   |\phi_{n+1-k,n}^{m}| 
&=& \sum_{m=1}^{\infty} |h_{n+m-1}| \sum_{k=1}^{n}  |\phi_{k,n}^{m}| \\
&\leq& \sum_{m=1}^{\infty} |h_{n+m-1}| \left\{ \sum_{k=1}^{n}  |\phi_{k}^{m}| +  \sum_{k=1}^{n}  |\phi_{k}^{m} - \phi_{k}^{m}|\right\}.
\end{eqnarray*} First, we assume (H-SM2). Then, 
\begin{equation*}
\sum_{m=1}^{\infty} \sum_{k=1}^{n}  |h_{n+m-1}|   |\phi_{n+1-k,n}^{m}| 
\leq \sup_{k \geq n} |h_k| \left\{ \sum_{m=1}^{\infty}  \sum_{k=1}^{n}  |\phi_{k}^{m}| +  \sum_{m=1}^{\infty} \sum_{k=1}^{n}  |\phi_{k}^{m} - \phi_{k}^{m}|\right\}.
\end{equation*} By (\ref{eq:unif-Baxter1}) and (\ref{eq:unif-Baxter3}), the second term above is $o(n^{-\alpha+1})$ as $n \rightarrow \infty$. Using similar expansion as in (\ref{eq:phikm-SM2}), for $n \in \N$, the first term is bounded by
\begin{equation*}
 \sum_{m=1}^{\infty}  \sum_{k=1}^{n}  |\phi_{k}^{m}| \leq  \sum_{m=1}^{\infty}  \sum_{k=1}^{n} \sum_{\ell=0}^{m-1} |\psi_\ell| |\phi_{k+m-1-\ell}| 
= \sum_{\ell=0}^{\infty} |\psi_\ell|  \sum_{k=1}^{n} \sum_{j=0}^{\infty} |\phi_{k+j}| 
\leq \|\psi\|_0 \|\phi\|_1 <\infty.
\end{equation*} All together, under (H-SM2),
\begin{eqnarray}
\sum_{m=1}^{\infty} \sum_{k=1}^{n}  |h_{n+m-1}|   |\phi_{n+1-k,n}^{m}|  
&\leq& \big\{ \|\psi\|_0 \|\phi\|_1  + o(n^{-\alpha+1}) \big\} \times \sup_{k \geq n} |h_k| \nonumber \\
&=& O(\sup_{k \geq n} |h_k|)
\quad \text{as}~~ n \rightarrow \infty.
 \label{eq:tail-SM2}
\end{eqnarray}

Next, we assume (H-SM1). Since $(h_1, h_2, ...) \in \ell_1(\N)$, we have
\begin{equation*}
 \sum_{m=1}^{\infty} |h_{n+m-1}| \sum_{k=1}^{n}  |\phi_{k,n}^{m}|
\leq \sum_{k=n}^{\infty} |h_{k}|
\left\{ \sup_{m \in \N} \sum_{k=1}^{n}  |\phi_{k}^{m}| + \sup_{m \in \N} \sum_{k=1}^{n}  |\phi_{k}^{m}-\phi_{k,n}^{m}|   \right\} 
.
\end{equation*} Using (\ref{eq:unif-Baxter1}), (\ref{eq:unif-Baxter2}), and (\ref{eq:phikm11}), it is easily seen that
\begin{equation*}
\sup_{m \in \N} \sum_{k=1}^{n}  |\phi_{k}^{m}| \leq 
 \|\psi\|_0 \|\phi\|_0
\quad \text{and} \quad 
\sup_{m \in \N} \sum_{k=1}^{n}  |\phi_{k}^{m}-\phi_{k,n}^{m}| = o(n^{-\alpha}) \quad \text{as}~~ n\rightarrow \infty.
\end{equation*} Therefore, we have
\begin{eqnarray}
 \sum_{m=1}^{\infty} |h_{n+m-1}| \sum_{k=1}^{n}  |\phi_{k,n}^{m}| 
&\leq& 
\left\{   \|\psi\|_0 \|\phi\|_0+ o(n^{-\alpha}) \right\} \times \sum_{k=n}^{\infty} |h_{k}| \nonumber \\
&=& O(\sum_{k=n}^{\infty} |h_{k}|)
 \quad \text{as}~~ n \rightarrow \infty.
\label{eq:tail-SM1}
\end{eqnarray}

Lastly, we assume (H-LM). Using that $(h_1, h_2, ...) \in \ell_1(\N)$, we have
\begin{equation*}
 \sum_{m=1}^{\infty} |h_{n+m-1}| \sum_{k=1}^{n}  |\phi_{k,n}^{m}|
\leq  \sum_{m=1}^{\infty} |h_{n+m-1}|\sum_{k=1}^{n}  |\phi_{k}^{m}-\phi_{k,n}^{m}|
+  \sum_{k=n}^{\infty} |h_{k}| \sup_{m \in \N} \sum_{k=1}^{n}  |\phi_{k}^{m}|.
\end{equation*} Using (\ref{eq:unif-BaxterLM1}), for large enough $n \in \N$, the first term above is bounded by
\begin{equation*}
\sum_{m=1}^{\infty} |h_{n+m-1}|\sum_{k=1}^{n}  |\phi_{k}^{m}-\phi_{k,n}^{m}|
\leq C_2 n^{-d} \sum_{m=1}^{\infty} |h_{n+m-1}| m^{d}
\leq C_2 n^{-d} \sum_{k=n}^{\infty} |h_{k}| k^{d}.
\end{equation*} Since $\{h_k\} \in \mathcal{C}_d$, the above is $C_2 \times o(n^{-d})$ as $n\rightarrow \infty$. 
Next, using (\ref{eq:unif-BaxterLM2}), we have
\begin{equation*}
\sum_{k=1}^{n}  |\phi_{k}^{m}| \leq \sum_{k=1}^{\infty}  |\phi_{k}^{m}| \leq C_3 \left(\frac{m}{m+1}\right)^d \leq C_3, \quad 
 m \in \N. 
\end{equation*}
Therefore, the second term is bounded by
\begin{equation*}
\sum_{k=n}^{\infty} |h_{k}| \sup_{m \in \N} |\sum_{k=1}^{n}  |\phi_{k}^{m}|
\leq C_3 n^{-d} \sum_{k=n}^{\infty} |h_{k}| k^{d} = C_3 \times o(n^{-d}) \quad \text{as}~~ n\rightarrow \infty.
\end{equation*} All together, we have
\begin{equation}\label{eq:tail-LM}
 \sum_{m=1}^{\infty} |h_{n+m-1}| \sum_{k=1}^{n}  |\phi_{k,n}^{m}| = \{C_2 + C_3\} \times o(n^{-d})
=o(n^{-d})
 \quad \text{as}~~ n \rightarrow \infty.
\end{equation}  Thus combining (\ref{eq:tail-SM2}), (\ref{eq:tail-SM1}), and (\ref{eq:tail-LM}) we get the desired results.
\hfill $\Box$

\begin{lemma} \label{lemma:abs-conv3}
Suppose $\{X_t\}$ satisfies either Assumption \ref{assum:SM} or \ref{assum:LM} and the filter coefficients $\{h_k\}$ satisfy Assumption \ref{assum:H}. 
Then,
\begin{equation} 
\bigg\| \sum_{k=n+1}^{\infty} \hat{h}_k X_{n+1-k} \bigg\|_V =
\begin{cases}
o(n^{-\alpha}) + O\big(\sum_{k=n}^{\infty} |h_k|\big), & \text{under (H-SM1)} \\
o(n^{-\alpha+1}) + O\big(\big\{ \sum_{k=n}^{\infty} |h_k|^2\big\} \big), & \text{under (H-SM2)} \\
O(n^{-d}), & \text{under (H-LM)} 
\end{cases} \quad \text{as}~~ n \rightarrow \infty.
\label{eq:Vbound}
\end{equation}
\end{lemma}
\noindent \textbf{Proof}. First, we assume (H-SM2). Since we assume short memory, we have
\begin{equation} \label{eq:spec-norm-ineq}
\bigg\| \sum_{k=n+1}^{\infty} \hat{h}_k X_{n+1-k} \bigg\|_V^2 \leq \sup_{\omega} f(e^{i\omega}) \sum_{k=n+1}^{\infty} |\hat{h}_k|^2,
\end{equation} where $f(\cdot)$ is a spectral density function of $\{X_t\}$. Using a similar expansion as in (\ref{eq:hsq11}), we have
\begin{eqnarray*} 
\sum_{k=n+1}^{\infty} |\hat{h}_k|^2 &=&
\sum_{k=n+1}^{\infty} \left\{ 
\sum_{m=1}^{\infty} |h_{-m}|  |\phi_{k}^m|  + |h_{k-1}|
\right\}^2 \\
&\leq&
 2  \sum_{m=1}^{\infty} |h_{-m}|^2 \cdot  \left\{ \sum_{k=n+1}^{\infty} \sum_{m=1}^{\infty} |\phi_{k}^m| \right\}^2
+ 2 \sum_{k=n+1}^{\infty} |h_{k-1}|^2.
\end{eqnarray*} Using (\ref{eq:phikm-SM2}), the above is bounded by
\begin{equation*}
\sum_{k=n+1}^{\infty} |\hat{h}_k|^2 = 2\|h\|_{\ell_2^-}^2 \|\psi\|_0^2 \times
o(n^{-2(\alpha-1)}) +2 \sum_{k=n}^{\infty} |h_{k}|^2 \quad
\text{as}~~ n\rightarrow \infty,
\end{equation*} where $\|h\|_{\ell_2^+}^2 =\sum_{m=1}^{\infty} |h_{-m}|^2$.
Substitute the above into (\ref{eq:spec-norm-ineq}) and using that
$\sup_{\omega} f(e^{i\omega}) \leq \sigma^2 \|\psi\|_0^2$ due to (\ref{eq:spec-decom1}), 
we have
\begin{eqnarray}
\bigg\| \sum_{k=n+1}^{\infty} \hat{h}_k X_{n+1-k} \bigg\|_V &\leq&
\sigma \|\psi\|_0 \bigg[
\sqrt{2} \|h\|_{\ell_2^-} \|\psi\|_0 \times
o(n^{-\alpha+1}) +\sqrt{2} \left\{\sum_{k=n}^{\infty} |h_{k}|^2\right\}^{1/2}\bigg] \nonumber \\
&=&
o(n^{-\alpha+1}) + O\left(\left\{\sum_{k=n}^{\infty} |h_{k}|^2\right\}^{1/2}\right) \quad \text{as}~~ n\rightarrow \infty.
 \label{eq:Vbound-SM2}
\end{eqnarray}

\vspace{0.5em}

Next, we assume (H-SM1) or (H-LM). By using (\ref{eq:FIRcoeff1}) together with a triangular inequality, we have
\begin{eqnarray}
\bigg\| \sum_{k=n+1}^{\infty} \hat{h}_k X_{n+1-k} \bigg\|_V &\leq& \sqrt{\var X_0} \sum_{k=n+1}^{\infty} |\hat{h}_k| \nonumber \\
&\leq& \sqrt{\var X_0} \left\{ \sum_{m=1}^{\infty} |h_{-m}|   \sum_{k=n+1}^{\infty} |\phi_{k}^m| + 
 \sum_{k=n+1}^{\infty} |h_{k-1}|\right\}.
 \label{eq:Vnorm-ineq}
\end{eqnarray} Under (H-SM1), by using (\ref{eq:phikm11}), the first term above is bounded by
\begin{equation*}
\sum_{m=1}^{\infty} |h_{-m}|   \sum_{k=n+1}^{\infty} |\phi_{k}^m| \leq
 \sup_{m \in \N}\sum_{k=n+1}^{\infty} |\phi_{k}^m| \cdot \sum_{m=1}^{\infty} |h_{-m}| 
=  \|h\|_{\ell_2^-} \|\psi\|_0 \times o(n^{-\alpha}) \quad
\text{as}~~ n\rightarrow \infty.
\end{equation*}Therefore, combining the above together with $\{\var X_0\}^{1/2} = \gamma(0)^{1/2} = \sigma \|\psi\|_{\ell_2^+}$,
where $\|\psi\|_{\ell_2^+} = \{ \sum_{j=0}^{\infty} |\psi_j|^2\}^{1/2}$,
we have
\begin{eqnarray}
\bigg\| \sum_{k=n+1}^{\infty} \hat{h}_k X_{n+1-k} \bigg\|_V &\leq& 
\sigma \|\psi\|_{\ell_2^+} \left[ 
\|h\|_{\ell_2^-} \|\psi\|_0 \times o(n^{-\alpha}) +  \sum_{k=n}^{\infty} |h_{k}|
\right] \nonumber \\
&=& o(n^{-\alpha+1}) + \sum_{k=n}^{\infty} |h_{k}|
\quad
\text{as}~~ n\rightarrow \infty.
 \label{eq:Vbound-SM1}
\end{eqnarray}
Under (H-LM), by using (\ref{eq:unif-BaxterLM2}), the first term in (\ref{eq:Vnorm-ineq}) is bounded by
\begin{eqnarray*}
\sum_{m=1}^{\infty} |h_{-m}|   \sum_{k=n+1}^{\infty} |\phi_{k}^m| \leq
C_3 \sum_{m=1}^{\infty} |h_{-m}| \left( \frac{m}{n+m}\right)^{d}
\leq C_3 n^{-d} \sum_{m=1}^{\infty} |h_{-m}| m^{d} = C_3 \times O(n^{-d})
\end{eqnarray*} as $n\rightarrow \infty$.
The second inequality is due to $(m/(n+m))^d \leq m^d/n^d$ for all $n,m \in \N$. Since, $\{h_k\} \in \mathcal{C}_d$, the second term is 
$ \sum_{k=n+1}^{\infty} |h_{k-1}| = o(n^{-d})$ as $n\rightarrow \infty$. All together,
\begin{equation} \label{eq:Vbound-LM}
\bigg\| \sum_{k=n+1}^{\infty} \hat{h}_k X_{n+1-k} \bigg\|_V = 
\sigma \|\psi\|_{\ell_2^+} \left[ C_3 \times O(n^{-d})  + o(n^{-d})\right]
=O(n^{-d})
 \quad
\text{as}~~ n\rightarrow \infty.
\end{equation} Combining (\ref{eq:Vbound-SM2}), (\ref{eq:Vbound-SM1}), and (\ref{eq:Vbound-LM}), we prove the lemma.
\hfill $\Box$

\end{appendix}


\bibliography{bib_LPP}
\bibliographystyle{plainnat}

\end{document}